# ON OPTIMAL SPATIAL SUBSAMPLE SIZE FOR VARIANCE ESTIMATION[1]

By Daniel J. Nordman and Soumendra N. Lahiri

*University of Wisconsin–La Crosse and Iowa State University*

We consider the problem of determining the optimal block (or subsample) size for a spatial subsampling method for spatial processes observed on regular grids. We derive expansions for the mean square error of the subsampling variance estimator, which yields an expression for the theoretically optimal block size. The optimal block size is shown to depend in an intricate way on the geometry of the spatial sampling region as well as characteristics of the underlying random field. Final expressions for the optimal block size make use of some nontrivial estimates of lattice point counts in shifts of convex sets. Optimal block sizes are computed for sampling regions of a number of commonly encountered shapes. Numerical studies are performed to compare subsampling methods as well as procedures for estimating the theoretically best block size.

**1. Introduction.** In this article, the problem of choosing subsample sizes is examined to maximize the performance of subsampling methods for variance estimation. The data at hand are viewed as realizations of a stationary, weakly dependent spatial lattice process. We consider the common scenario of sampling from sites of regular distance (e.g., indexed by the integer lattice $\mathbb{Z}^d$), lying within some region $R_n$ embedded in $\mathbb{R}^d$. Such lattice data appear often in time series, agricultural field trials, and remote sensing and image analysis (medical and satellite image processing).

Consider estimating the variance of a statistic $\hat{\theta}_n$ from $R_n$. For variance estimation via subsampling, the basic idea is to construct several "scaled-down" copies (subsamples) of the sampling region $R_n$ that fit inside $R_n$, evaluate the analog of $\hat{\theta}_n$ on each of these subregions, and then compute

Received July 2002; revised July 2003.
[1]Supported in part by NSF Grants DMS-98-16630 and DMS-00-72571 and by the Deutsche Forschungsgemeinschaft (SFB 475).

a properly normalized sample variance from the resulting values. The $R_n$-sampling scheme is essentially recreated at the level of the subregions. Two subsampling designs are most typical: Subregions can be maximally overlapping (OL) or devised to be nonoverlapping (NOL). The accuracy (e.g., variance and bias) of subsample-based estimators depends crucially on the choice of subsample size.

To place our work into perspective, we briefly outline previous research in variance estimation with subsamples and theoretical size considerations. Variance estimation through subsampling originated from analysis of weakly dependent, stationary time processes. Suppose $\hat{\theta}_n$ is an estimator of a parameter of interest $\theta$ based on $\{Z(1), \ldots, Z(n)\}$ from a stationary temporal process $\{Z(i)\}_{i \geq 1}$. To obtain subsamples for $\hat{\theta}_n$-variance estimation, Carlstein (1986) first proposed the use of NOL blocks of length $m \leq n$: $\{Z(1+(i-1)m), \ldots, Z(im)\}$, $i = 1, \ldots, \lfloor n/m \rfloor$, while the sequence of subseries $\{Z(i), \ldots, Z(i+m-1)\}$, $i = 1, \ldots, n-m+1$, provides OL subsamples of length $m$ [cf. Künsch (1989) and Politis and Romano (1993b)]. Here, $\lfloor x \rfloor$ denotes the integer part of a real number $x$. In each respective subsample collection, evaluations of an analog statistic $\hat{\theta}_i$ are made for each subseries and a normalized sample variance is calculated to estimate the parameter $n \operatorname{Var}(\hat{\theta}_n)$,

$$\sum_{i=1}^{J} \frac{m(\hat{\theta}_i - \tilde{\theta})^2}{J}, \qquad \tilde{\theta} = \sum_{i=1}^{J} \frac{\hat{\theta}_i}{J},$$

where $J = \lfloor n/m \rfloor$ ($J = n - m + 1$) for the NOL (OL) subsample-based estimator. Carlstein (1986) and Fukuchi (1999) established the $L_2$ consistency of the NOL and OL estimators, respectively, for the variance of a general (not necessarily linear) statistic. Politis and Romano (1993b) determined asymptotic orders of the variance $O(m/n)$ and bias $O(1/m)$ of the subsample variance estimators for linear statistics. For mixing time series, they found that a subsample size $m$ proportional to $n^{1/3}$ is optimal in the sense of minimizing the mean square error (MSE) of variance estimation, concurring also with optimal block order for the moving block bootstrap variance estimator [Hall, Horowitz and Jing (1995) and Lahiri (1996)].

Cressie [(1991), page 492] conjectured the recipe for extending Carlstein's variance estimator to the general spatial setting, obtaining subsamples by tiling the sample region $R_n$ with disjoint "congruent" subregions. Politis and Romano (1993a, 1994) have shown the consistency of subsample-based variance estimators for rectangular sampling or subsampling regions in $\mathbb{R}^d$ when the sampling sites are observed on $\mathbb{Z}^d \cap \prod_{i=1}^{d}[1, n_i]$ and integer translates of $\prod_{i=1}^{d}[1, m_i]$ yield the subsamples. Garcia-Soidan and Hall (1997) and Possolo (1991) proposed similar estimators under an identical sampling scenario. For linear statistics, Politis and Romano (1993a) determined that a subsampling scaling choice

$$\prod_{i=1}^{d} m_i = C \left\{ \prod_{i=1}^{d} n_i \right\}^{d/(d+2)},$$

for some unknown $C$, minimizes the order of a variance estimator's asymptotic MSE. Sherman and Carlstein (1994) and Sherman (1996) proved the MSE-consistency of NOL and OL subsample estimators, respectively, for the variance of general statistics in $\mathbb{R}^2$. Their work allowed for a more flexible sampling scheme: the "inside" of a simple closed curve defines a set $D \subset [-1,1]^2$, $\mathbb{Z}^2 \cap nD$ (using a scaled-up copy of $D$) constitutes the set of sampling sites, and translates of $mD$ within $nD$ form subsamples. Sherman (1996) minimized a bound on the asymptotic order of the OL estimator's MSE to argue that the best size choice for OL subsamples involves $m = O(n^{1/2})$ [coinciding with the above findings of Politis and Romano (1993a) for rectangular regions in $\mathbb{R}^2$]. Politis and Sherman (2001) have developed consistent subsampling methods for variance estimation with marked point process data [cf. Politis, Romano and Wolf (1999), Chapter 6].

Few theoretical and numerical recommendations for choosing subsamples have been offered in the spatial setting, especially with the intent of variance estimation. As suggested in the literature, an explicit theoretical determination of optimal subsample size or scaling requires calculation of an order and associated proportionality constant for a given sampling region $R_n$. Even for the few sampling situations where the order of optimal subsample size has been established, the exact adjustments to these orders are unknown and, quoting Politis and Romano (1993a), "important (and difficult) in practice." Beyond the time series case with the univariate sample mean, the influence of the geometry and dimension of $R_n$, as well as the structure of $\hat{\theta}_n$, on precise subsample selection has not been explored. We attempt here to advance some ideas on the best size choice, both theoretically and empirically, for subsamples.

We work under the "smooth function" model of Hall (1992), where the statistic of interest $\hat{\theta}_n$ can be represented as a function of sample means. We formulate a framework for sampling in $\mathbb{R}^d$ where the sampling region $R_n$ is obtained by "inflating" a prototype set in the unit cube in $\mathbb{R}^d$ and the subsampling regions are given by suitable translates of a scaled down copy of the sampling region $R_n$. We consider both a nonoverlapping version and a (maximal) overlapping version of the subsampling method. For each method, we derive expansions for the variance and the bias of the corresponding subsample estimator of $\text{Var}(\hat{\theta}_n)$. The asymptotic variance of the spatial subsample estimator for the OL version turns out to be smaller than that of the NOL version by a constant factor $K_1$ (say) which depends solely on the geometry of the sampling region $R_n$. In the time series case, Meketon and Schmeiser (1984), Künsch (1989), Hall, Horowitz and Jing (1995) and Lahiri (1996) have shown in different degrees of generality that the asymptotic variance under the OL subsampling scheme, compared to the NOL one, is $K_1 = \frac{2}{3}$ times smaller. Results of this paper show that for rectangular sampling regions $R_n$ in $d$-dimensional space, the factor $K_1$ is given by $(\frac{2}{3})^d$. We list the factor $K_1$ for sampling regions of some common shapes in Table 1.

In contrast, the bias parts of both the OL and NOL subsample variance estimators are usually asymptotically equivalent and depend on the covariance structure of the random field as well as on the geometry of the sampling region $R_n$. Since the bias term is typically of the same order as the number of lattice points lying near a subsample's boundary, determination of the leading bias term involves some nontrivial estimates of the lattice point counts over translated subregions. Counting lattice points in scaled-up sets is a hard problem and has received a lot of attention in analytic number theory and in combinatorics. Even for the case of the plane (i.e., $d = 2$), the counting results available in the literature are directly applicable to our problem only for a very restricted class of subregions that have the so-called "smoothly winding border" [cf. van der Corput (1920) and Huxley (1993, 1996)]. Here explicit expressions for the bias terms are derived for a more general class of sampling regions using some new estimates on the discrepancy between the number of lattice points and the volume of the *shifted* subregions in the plane and in three-dimensional Euclidean space. In particular, our results are applicable to sampling regions that do not necessarily have "smoothly winding borders."

Minimizing the combined expansions for the bias and the variance parts, we derive explicit expressions for the theoretical optimal block size for sampling regions of different shapes. To briefly describe the result for a few common shapes: Suppose the sampling region $R_n$ is obtained by inflating a given set $R_0 \in (-\frac{1}{2}, \frac{1}{2}]^d$ by a scaling constant $\lambda_n$ as $R_n = \lambda_n R_0$ and that the subsamples are formed by considering the translates of $_sR_n = {}_s\lambda_n R_0$. Then the theoretically optimal choice of the subsample size $_s\lambda_n$ for the OL version is of the form

$$_s\lambda_n^{\text{opt}} = \left(\frac{\lambda_n^d B_0^2}{dK_0\tau^4}\right)^{1/(d+2)} (1 + o(1)) \qquad \text{as } n \to \infty$$

Table 1
*Examples of $K_1$ for several shapes of the sampling region $R_n \subset \mathbb{R}^d$*

| Shape of $R_n$ | Rectangle in $\mathbb{R}^d$ | Sphere in $\mathbb{R}^3$ | Circle in $\mathbb{R}^2$ | Right triangle in $\mathbb{R}^2$ |
|---|---|---|---|---|
| $K_1$ | $(2/3)^d$ | $17\pi/315$ | $\pi/4 - 4/(3\pi)$ | $1/5$ |

for some constants $B_0$ and $K_0$ (coming from the bias and the variance terms, respectively) where $\tau^2$ is a population parameter that does not depend on the shape of the sampling region $R_n$ (see Theorem 5.1 for details). Table 2 lists the constants $B_0$ and $K_0$ for some shapes of $R_n$. It follows from Table 2 that, unlike the time series case, in higher dimensions the optimal block size critically depends on the shape of the spatial sampling region $R_n$. It simplifies only slightly for the NOL subsampling scheme as the constant $K_0$ is unnecessary for computing optimal NOL subsamples, but the bias constant $B_0$ is often the same for estimators from each version of subsampling. These expressions may be readily used to obtain estimates of the theoretical optimal subsample scaling for use in practice.

The rest of the paper is organized as follows. In Section 2 we describe the spatial subsampling method and state the assumptions used in the paper. In Sections 3 and 4 we, respectively, derive expansions for the variance and the bias parts of the subsampling estimators. Theoretical optimal subsample scalings (or block sizes) are derived in Section 5. The results are illustrated with some common examples in Section 6. Section 7 describes two methods for estimating optimal subsample scaling. In Section 8 a numerical study of subsample variance estimators and scaling estimation methods is provided. Proofs of variance and bias results are separated into Sections 9 and 10, respectively.

**2. Variance estimators via subsampling.** In Section 2.1 we frame the sampling design and the structure of the sampling region. Two methods of subsampling are presented in Section 2.2 along with corresponding nonparametric variance estimators. Assumptions and conditions used in the paper are given in Section 2.3.

2.1. *The sampling structure.* To describe the sampling scheme used, we first assume all potential sampling sites are located on a translate of the rectangular integer lattice in $\mathbb{R}^d$. For a fixed (chosen) vector $\mathbf{t} \in [-1/2, 1/2]^d$, we identify the $\mathbf{t}$-translated integer lattice as $\mathbf{Z}^d \equiv \mathbf{t} + \mathbb{Z}^d$. Let $\{Z(\mathbf{s}) : \mathbf{s} \in \mathbf{Z}^d\}$ be a stationary weakly dependent random field (hereafter r.f.) taking values in $\mathbb{R}^p$. [We use bold font as a standard to denote vectors in the space of sampling $\mathbb{R}^d$ and normal font for vectors in $\mathbb{R}^p$, including $Z(\cdot)$.] We suppose that the process $Z(\cdot)$ is observed at sampling sites lying within the sampling region $R_n \subset \mathbb{R}^d$. That is, the collection of available sampling sites is $\{Z(\mathbf{s}) : \mathbf{s} \in R_n \cap \mathbf{Z}^d\}$.

TABLE 2
*Examples of $B_0$, $K_0$ for some sampling regions $R_n$*[*]

| $R_n$ | Sphere in $\mathbb{R}^3$ | Cross in $\mathbb{R}^2$  | Right triangle in $\mathbb{R}^2$  |
|---|---|---|---|
| $K_0$ | $34/105$ | $4/9 \cdot 191/192$ | $2/5$ |
| $B_0$ | $3/2 \sum_{\mathbf{k}\in\mathbb{Z}^3} \|\mathbf{k}\|\sigma(\mathbf{k})$ | $4/3 \sum_{\mathbf{k}\in\mathbb{Z}^2} \|\mathbf{k}\|_1 \sigma(\mathbf{k})$ | $2\sum_{\mathbf{k}=(k_1,k_2)'\in\mathbb{Z}^2, \operatorname{sign} k_1 = \operatorname{sign} k_2} \|\mathbf{k}\|_1 \sigma(\mathbf{k})$ $+ 2\sum_{\mathbf{k}\in\mathbb{Z}^2, \operatorname{sign} k_1 \neq \operatorname{sign} k_2} \|\mathbf{k}\|_\infty \sigma(\mathbf{k})$ |

[*]Cross and triangle shapes appear in Figure 1; see Section 6 for further details. Autocovariances $\sigma(\cdot)$ and Euclidean, $l^1$, and $l^\infty$ norms $\|\cdot\|$, $\|\cdot\|_1$, $\|\cdot\|_\infty$ are described in Section 2.3.

To obtain the results in the paper, we assume that the sampling region $R_n$ becomes unbounded as the sample size increases. This will provide a commonly used "increasing domain" framework for studying asymptotics with spatial lattice data [cf. Cressie (1991)]. We next specify the structure of the regions $R_n$ and employ a formulation similar to that of Lahiri (1999a, 2004).

Let $R_0$ be a Borel subset of $(-1/2, 1/2]^d$ containing an open neighborhood of the origin such that for any sequence of positive real numbers $a_n \to 0$, the number of cubes of the scaled lattice $a_n \mathbb{Z}^d$ which intersect the closures $\overline{R_0}$ and $\overline{R_0^c}$ is $O((a_n^{-1})^{d-1})$ as $n \to \infty$. Let $\Delta_n$ be a sequence of $d \times d$ diagonal matrices, with positive diagonal elements $\lambda_1^{(n)}, \ldots, \lambda_d^{(n)}$, such that each $\lambda_i^{(n)} \to \infty$ as $n \to \infty$. We assume that the sampling region $R_n$ is obtained by "inflating" the template set $R_0$ by the directional scaling factors $\Delta_n$; namely,

$$R_n = \Delta_n R_0.$$

Because the origin is assumed to lie in $R_0$, the sampling region $R_n$ grows outward in all directions as $n$ increases. Furthermore, if the scaling factors are all equal $(\lambda_1^{(n)} = \cdots = \lambda_d^{(n)})$, the shape of $R_n$ remains the same for different values of $n$.

The formulation given above allows the sampling region $R_n$ to have a large variety of fairly irregular shapes with the boundary condition on $R_0$ imposed to avoid pathological cases. Some common examples of such regions are convex subsets of $\mathbb{R}^d$, such as spheres, ellipsoids, polyhedrons, as well as certain nonconvex subsets with irregular boundaries, such as star-shaped regions. Sherman and Carlstein (1994) and Sherman (1996) consider a similar class of such regions in the plane (i.e., $d = 2$) where the boundaries of the sets $R_0$ are delineated by simple rectifiable curves with finite lengths. The border requirements on $R_0$ ensure that the number of observations near the boundary of $R_n$ is negligible compared to the totality of data values.

2.2. *Subsampling designs and variance estimators.* We suppose that the relevant statistic, whose variance we wish to estimate, can be represented as a function of sample means. Let $\hat{\theta}_n = H(\bar{Z}_{N_n})$ be an estimator of the population parameter of interest $\theta = H(\mu)$, where $H:\mathbb{R}^p \to \mathbb{R}$ is a smooth function, $EZ(\mathbf{t}) = \mu \in \mathbb{R}^p$ is the mean of the stationary r.f., and $\bar{Z}_{N_n}$ is the sample mean of the $N_n$ observations within $R_n$,

$$\bar{Z}_{N_n} = N_n^{-1} \sum_{\mathbf{s} \in \mathbf{Z}^d \cap R_n} Z(\mathbf{s}). \tag{2.1}$$

This parameter and estimator formulation is what Hall (1992) calls the "smooth function" model and it has been used in other scenarios, such as with the moving block bootstrap (MBB), for studying approximately linear functions of a sample mean [cf. Lahiri (1996) and Politis, Romano and Wolf (1999)]. By considering suitable functions of the $Z(\mathbf{s})$'s, one can represent a wide range of estimators under the present framework. In particular, these include means, products and ratios of means, sample moments, spatial correlograms, Yule–Walker estimates for autoregressive processes [cf. Guyon (1995)] and some pseudo likelihood-based estimators of process parameters [cf. Ripley (1981)].

The quantity which we seek to estimate nonparametrically is the variance of the normalized statistic $\sqrt{N_n}\hat{\theta}_n$, say, $\tau_n^2 = N_n E(\hat{\theta}_n - E\hat{\theta}_n)^2$. In our problem, this goal is equivalent to consistently estimating the limiting variance $\tau^2 = \lim_{n \to \infty} \tau_n^2$.

2.2.1. *Overlapping subsamples.* Variance estimation with OL subsampling regions has often been considered in the literature, though in more narrow sampling situations [cf. $\mathbb{R}^2$-sampling regions, Sherman (1996); $\mathbb{R}^d$-rectangular regions, Politis and Romano (1994); time series data, Politis and Romano (1993a)].

We first consider creating a smaller version of $R_n$, which will serve as a template for the OL subsampling regions. To this end, let $_s\Delta_n$ be a $d \times d$ diagonal matrix with positive diagonal elements, $\{_s\lambda_1^{(n)}, \ldots, _s\lambda_d^{(n)}\}$, such that $_s\lambda_i^{(n)}/\lambda_i^{(n)} \to 0$ and $_s\lambda_i^{(n)} \to \infty$, as $n \to \infty$, for each $i = 1, \ldots, d$. (The matrix $\Delta_n$ represents the determining scaling factors for $R_n$ and $_s\Delta_n$ shall be factors used to define the subsamples.) We make the "prototype" subsampling region

$$_sR_n = {_s\Delta_n} R_0, \tag{2.2}$$

and identify a subset of $\mathbb{Z}^d$, say $J_{\text{OL}}$, corresponding to all integer translates of $_sR_n$ lying within $R_n$. That is,

$$J_{\text{OL}} = \{\mathbf{i} \in \mathbb{Z}^d : \mathbf{i} + {_sR_n} \subset R_n\}.$$

The desired OL subsampling regions are precisely the translates of $_sR_n$ given by $R_{\mathbf{i},n} \equiv \mathbf{i} + {_sR_n}$, $\mathbf{i} \in J_{\text{OL}}$. Note that the origin belongs to $J_{\text{OL}}$ and some of these subregions may clearly overlap.

Let $_sN_n = |\mathbf{Z}^d \cap {_sR_n}|$ be the number of sampling sites in $_sR_n$ and let $|J_{\mathrm{OL}}|$ denote the number of available subsampling regions. The number of sampling sites within each OL subsampling region is the same, namely for any $\mathbf{i} \in J_{\mathrm{OL}}$, $_sN_n = |\mathbf{Z}^d \cap R_{\mathbf{i},n}|$. For each $\mathbf{i} \in J_{\mathrm{OL}}$, compute $\hat{\theta}_{\mathbf{i}}^{\mathrm{OL}} = H(Z_{\mathbf{i},n})$, where

$$Z_{\mathbf{i},n} = {_sN_n}^{-1} \sum_{\mathbf{s} \in \mathbf{Z}^d \cap R_{\mathbf{i},n}} Z(\mathbf{s})$$

denotes the sample mean of observations within the subregion. We then have the OL subsample variance estimator of $\tau_n^2$ as

$$\hat{\tau}_{n,\mathrm{OL}}^2 = |J_{\mathrm{OL}}|^{-1} \sum_{\mathbf{i} \in J_{\mathrm{OL}}} {_sN_n}(\hat{\theta}_{\mathbf{i},n}^{\mathrm{OL}} - \tilde{\theta}_n^{\mathrm{OL}})^2,$$

$$\tilde{\theta}_n^{\mathrm{OL}} = |J_{\mathrm{OL}}|^{-1} \sum_{\mathbf{i} \in J_{\mathrm{OL}}} \hat{\theta}_{\mathbf{i},n}^{\mathrm{OL}}.$$

2.2.2. *Nonoverlapping subsamples.* To create NOL subsamples, we adopt a formulation similar to that of Sherman and Carlstein (1994) and Lahiri (1999a). The sampling region $R_n$ is first divided into disjoint "cubes." Let $_s\Delta_n$ be the previously described $d \times d$ diagonal matrix from (2.2), which will determine the "window width" of the partitioning cubes. Let

$$J_{\mathrm{NOL}} = \{\mathbf{i} \in \mathbb{Z}^d : {_s\Delta_n}(\mathbf{i} + (-1/2, 1/2]^d) \subset R_n\}$$

represent the set of all "inflated" subcubes that lie inside $R_n$. Denote its cardinality as $|J_{\mathrm{NOL}}|$. For each $\mathbf{i} \in J_{\mathrm{NOL}}$, define the subsampling region $\tilde{R}_{\mathbf{i},n} = {_s\Delta_n}(\mathbf{i} + R_0)$ by inscribing the translate of $_s\Delta_n R_0$ such that the origin is mapped onto the midpoint of the cube $_s\Delta_n(\mathbf{i} + (-1/2, 1/2]^d)$. This provides a collection of NOL subsampling regions, which are smaller versions of the original sampling region $R_n$ that lie inside $R_n$.

For each $\mathbf{i} \in J_{\mathrm{NOL}}$, the function $H(\cdot)$ is evaluated at the sample mean, say $\tilde{Z}_{\mathbf{i},n}$, for a corresponding subsampling region $\tilde{R}_{\mathbf{i},n}$ to obtain $\hat{\theta}_{\mathbf{i},n}^{\mathrm{NOL}} = H(\tilde{Z}_{\mathbf{i},n})$. The NOL subsample estimator of $\tau_n^2$ is again an appropriately scaled sample variance,

$$\hat{\tau}_{n,\mathrm{NOL}}^2 = |J_{\mathrm{NOL}}|^{-1} \sum_{\mathbf{i} \in J_{\mathrm{NOL}}} {_sN_{\mathbf{i},n}}(\hat{\theta}_{\mathbf{i},n}^{\mathrm{NOL}} - \tilde{\theta}_n^{\mathrm{NOL}})^2,$$

$$\tilde{\theta}_n^{\mathrm{NOL}} = |J_{\mathrm{NOL}}|^{-1} \sum_{\mathbf{i} \in J_{\mathrm{NOL}}} \hat{\theta}_{\mathbf{i},n}^{\mathrm{NOL}},$$

where $_sN_{\mathbf{i},n} = |\mathbf{Z}^d \cap \tilde{R}_{\mathbf{i},n}|$ denotes the number of sampling sites within a given NOL subsample.

We note that $_sN_{\mathbf{i},n}$ may differ between NOL subsamples, but all such subsamples will have exactly $_sN_{\mathbf{i},n} = {_sN_n}$ sites available if the diagonal elements of $_s\Delta_n$ are integers.

2.3. *Assumptions.* For stating the assumptions, we need to introduce some notation. For a vector $\mathbf{x} = (x_1, \ldots, x_d)' \in \mathbb{R}^d$, let $\|\mathbf{x}\|$ and $\|\mathbf{x}\|_1 = \sum_{i=1}^d |x_i|$ denote the usual Euclidean and $l^1$ norms of $\mathbf{x}$, respectively. Denote the $l^\infty$ norm as $\|\mathbf{x}\|_\infty = \max_{1 \leq k \leq d} |x_k|$. Define $\operatorname{dis}(E_1, E_2) = \inf\{\|\mathbf{x} - \mathbf{y}\|_\infty : \mathbf{x} \in E_1, \mathbf{y} \in E_2\}$ for two sets $E_1, E_2 \subset \mathbb{R}^d$. We shall use the notation $|\cdot|$ also in two other cases: for a countable set $B$, $|B|$ will denote the cardinality of the set $B$; for an uncountable set $A \subset \mathbb{R}^d$, $|A|$ will refer to the volume (i.e., the $\mathbb{R}^d$ Lebesgue measure) of $A$.

Let $\mathcal{F}_Z(T) = \sigma\langle Z(\mathbf{s}) : \mathbf{s} \in T\rangle$ be the $\sigma$-field generated by the variables $\{Z(\mathbf{s}) : \mathbf{s} \in T\}$, $T \subset \mathbf{Z}^d$. For $T_1, T_2 \subset \mathbf{Z}^d$, write $\tilde{\alpha}(T_1, T_2) = \sup\{|P(A \cap B) - P(A)P(B)| : A \in \mathcal{F}_Z(T_1), B \in \mathcal{F}_Z(T_2)\}$. Then the strong mixing coefficient for the r.f. $Z(\cdot)$ is defined as

(2.3) $\alpha(k, l) = \sup\{\tilde{\alpha}(T_1, T_2) : T_i \subset \mathbf{Z}^d,\ |T_i| \leq l,\ i = 1, 2;\ \operatorname{dis}(T_1, T_2) \geq k\}.$

Note that the supremum in the definition of $\alpha(k, l)$ is taken over sets $T_1, T_2$ which are bounded. For $d > 1$ this is important. An r.f. on the lattice $\mathbb{Z}^d$ with $d \geq 2$ that satisfies a strong mixing condition of the form

(2.4) $\lim_{k \to \infty} \sup\{\tilde{\alpha}(T_1, T_2) : T_1, T_2 \subset \mathbf{Z}^d,\ \operatorname{dis}(T_1, T_2) \geq k\} = 0$

with supremum taken over possibly unbounded sets necessarily belongs to the more restricted class of $\rho$-mixing r.f.'s [cf. Bradley (1989)]. Politis and Romano (1993a) use moment inequalities based on the mixing condition in (2.4) to determine the orders of the bias and variance of $\hat{\tau}_{n,\mathrm{OL}}^2$, $\hat{\tau}_{n,\mathrm{NOL}}^2$ for rectangular sampling regions.

For proving the subsequent theorems, Assumptions A.1–A.5 are needed along with two conditions stated as functions of a positive argument $r \in \mathbb{Z}_+ = \{0, 1, 2, \ldots\}$. In the following, $\det(\Delta)$ represents the determinant of a square matrix $\Delta$. For $\alpha = (\alpha_1, \ldots, \alpha_p)' \in (\mathbb{Z}_+)^p$, let $D^\alpha$ denote the $\alpha$th order partial differential operator $\partial^{\alpha_1 + \cdots + \alpha_p}/\partial x_1^{\alpha_1} \cdots \partial x_p^{\alpha_p}$ and $\nabla = (\partial H(\mu)/\partial x_1, \ldots, \partial H(\mu)/\partial x_p)'$ be the vector of first-order partial derivatives of $H$ at $\mu$. Limits in order symbols are taken letting $n$ tend to infinity.

ASSUMPTIONS.

A.1. There exists a $d \times d$ diagonal matrix $\Delta_0$, $\det(\Delta_0) > 0$, such that

$$\frac{1}{{}_s\lambda_1^{(n)}} {}_s\Delta_n \to \Delta_0.$$

A.2. For the scaling factors of the sampling and subsampling regions

$$\sum_{i=1}^d \frac{1}{{}_s\lambda_i^{(n)}} + \sum_{i=1}^d \frac{{}_s\lambda_i^{(n)}}{\lambda_i^{(n)}} + \frac{[\det({}_s\Delta_n)]^{(d+1)/d}}{\det(\Delta_n)} = o(1),$$

$$\max_{1 \leq i \leq d} \lambda_i^{(n)} = O\Big(\min_{1 \leq i \leq d} \lambda_i^{(n)}\Big).$$

A.3. There exist nonnegative functions $\alpha_1(\cdot)$ and $g(\cdot)$ such that $\lim_{k\to\infty}\alpha_1(k) = 0$, $\lim_{l\to\infty}g(l) = \infty$ and the strong-mixing coefficient $\alpha(k,l)$ from (2.3) satisfies the inequality

$$\alpha(k,l) \leq \alpha_1(k)g(l), \qquad k > 0, \, l > 0.$$

A.4. $\sup\{\tilde{\alpha}(T_1,T_2) : T_1, T_2 \subset \mathbf{Z}^d, \, |T_1| = 1, \, \mathrm{dis}(T_1,T_2) \geq k\} = o(k^{-d})$.

A.5. $\tau^2 > 0$, where $\tau^2 = \sum_{\mathbf{k}\in\mathbb{Z}^d}\sigma(\mathbf{k})$, $\sigma(\mathbf{k}) = \mathrm{Cov}(\nabla' Z(\mathbf{t}), \nabla' Z(\mathbf{t}+\mathbf{k}))$.

CONDITIONS.

$D_r$. $H : \mathbb{R}^p \to \mathbb{R}$ is $r$-times continuously differentiable and, for some $a \in \mathbb{Z}_+$ and real $\mathcal{C} > 0$,

$$\max\{|D^\nu H(\mathbf{x})| : \|\nu\|_1 = r\} \leq \mathcal{C}(1 + \|\mathbf{x}\|^a), \qquad \mathbf{x} \in \mathbb{R}^p.$$

$M_r$. For some $0 < \delta \leq 1$, $0 < \kappa < (2r - 1 - 1/d)(2r + \delta)/\delta$, and $\mathcal{C} > 0$,

$$\mathrm{E}\|Z(\mathbf{t})\|^{2r+\delta} < \infty,$$

$$\sum_{m=1}^{\infty} m^{(2r-1)d-1}\alpha_1(m)^{\delta/(2r+\delta)} < \infty,$$

$$g(x) \leq \mathcal{C}x^\kappa, \qquad x \in [1, \infty).$$

Some comments about the assumptions and the conditions are in order. Assumption A.5 implies a positive, finite asymptotic variance $\tau^2$ for the standardized estimator $\sqrt{N_n}\hat{\theta}_n$.

In Assumption A.3 we formulate a conventional bound on the mixing coefficient $\alpha(k,l)$ from (2.3) that is applicable to many r.f.'s and resembles the mixing assumption of Lahiri (1999a, 2004). For r.f.'s satisfying Assumption A.3, the "distance" component of the bound, $\alpha_1(\cdot)$, often decreases at an exponential rate while the function of "set size," $g(\cdot)$, increases at a polynomial rate [cf. Guyon (1995)]. Examples of r.f.'s that meet the requirements of Assumption A.3 and Condition $M_r$ include Gaussian fields with analytic spectral densities, certain linear fields with a moving average or autoregressive (AR) representation (like $m$-dependent fields), separable $\mathrm{AR}(1) \times \mathrm{AR}(1)$ lattice processes suggested by Martin (1990) for modeling in $\mathbb{R}^2$, many Gibbs and Markov fields, and important time series models [cf. Doukhan (1995)]. Condition $M_r$ combined with Assumption A.3 also provides useful moment bounds for normed sums of observations (see Lemma 9.2).

Assumption A.4 permits the CLT in Bolthausen (1982) to be applied to sums of $Z(\cdot)$ on sets of increasing domain, in conjunction with the boundary condition on $R_0$, Assumption A.3 and Condition $M_r$. This version of the CLT (Stein's method) is derived from $\alpha$-mixing conditions which ensure asymptotic independence between a single point and observations in arbitrary sets of increasing distance [cf. Perera (1997)].

Assumptions A.1 and A.2 set additional guidelines for how sampling and subsampling design parameters, $\Delta_n$ and $_s\Delta_n$, may be chosen. The assumptions provide a flexible framework for handling "increasing domains" of many shapes. For $d=1$, Assumptions A.1 and A.2 are equivalent to the requirements of Lahiri (1999b) who provides variance and bias expansions for the MBB variance estimator with weakly dependent time processes.

**3. Variance expansions.** We now give expansions for the asymptotic variance of the OL/NOL subsample variance estimators $\hat{\tau}^2_{n,\text{OL}}$ and $\hat{\tau}^2_{n,\text{NOL}}$ of $\tau^2_n = N_n \text{Var}(\hat{\theta}_n)$.

THEOREM 3.1. *Suppose that Assumptions* A.1–A.5 *and Conditions* $D_2$ *and* $M_{5+2a}$ *hold with a as specified under Condition* $D_2$. *Then,*

(a) $$\text{Var}(\hat{\tau}^2_{n,\text{OL}}) = K_0 \cdot \frac{\det(_s\Delta_n)}{\det(\Delta_n)}[2\tau^4](1+o(1)),$$

(b) $$\text{Var}(\hat{\tau}^2_{n,\text{NOL}}) = \frac{1}{|R_0|} \cdot \frac{\det(_s\Delta_n)}{\det(\Delta_n)}[2\tau^4](1+o(1)),$$

*where*

$$K_0 = \frac{1}{|R_0|} \cdot \int_{\mathbb{R}^d} \frac{|(\mathbf{x}+R_0) \cap R_0|^2}{|R_0|^2} \, d\mathbf{x}$$

*is an integral with respect to the* $\mathbb{R}^d$ *Lebesgue measure.*

The constant $K_0$ appearing in the variance expansion of the estimator $\hat{\tau}^2_{n,\text{OL}}$ is a property of the *shape* of the sampling template $R_0$ but not of its exact embedding in space $\mathbb{R}^d$ or even the scale of the set. Namely, $K_0$ is invariant to invertible affine transformations applied to $R_0$ and hence can be computed from either $R_0$ or $R_n = \Delta_n R_0$. Values of $K_0$ for some template shapes are given in Table 3 and Section 6.

A stationary time sequence $Z(1), \ldots, Z(n)$ can be obtained within our sampling formulation by choosing $R_0 = (-1/2, 1/2]$ and $\lambda_1^{(n)} = n$ on the untranslated integer lattice $\mathbf{Z} = \mathbb{Z}$. In this special sampling case, an application

TABLE 3
*Examples of $K_0$ from Theorem* 3.1 *for several shapes of $R_0 \subset \mathbb{R}^d$*

| $R_0$ Shape | $\mathbb{R}^d$ Rectangle | $\mathbb{R}^3$ Ellipsoid | $\mathbb{R}^3$ Cylinder | $\mathbb{R}^2$ Ellipse | $\mathbb{R}^2$ Trapezoid* |
|---|---|---|---|---|---|
| $K_0$ | $(2/3)^d$ | $34/105$ | $2/3(1 - 16/(3\pi^2))$ | $1 - 16/(3\pi^2)$ | $2/5(1+4c/9)$ |

*The trapezoid has a 90° interior $\angle$ and parallel sides $b_2 \geq b_1$; $c = (b_2/b_1+1)^{-2}[1+2(b_2/b_1-1)/(b_2/b_1+1)]$.

of Theorem 3.1 yields

$$\mathrm{Var}(\hat{\tau}_{n,\mathrm{OL}}^2) = 2/3 \cdot \mathrm{Var}(\hat{\tau}_{n,\mathrm{NOL}}^2),$$

$$\mathrm{Var}(\hat{\tau}_{n,\mathrm{NOL}}^2) = {}_s\lambda_1^{(n)} \cdot [2\tau^4](1+o(1)),$$

a result which is well known for "nearly" linear functions $\hat{\theta}_n$ of a time series sample mean [cf. Künsch (1989)]. Theorem 3.1 implies that, under the "smooth" function model, the asymptotic variance of the OL subsample-based variance estimator is always strictly less than the NOL version because

$$(3.1) \qquad K_1 = \lim_{n \to \infty} \frac{\mathrm{Var}(\hat{\tau}_{n,\mathrm{OL}}^2)}{\mathrm{Var}(\hat{\tau}_{n,\mathrm{NOL}}^2)} = K_0 |R_0| < 1.$$

If both estimators have the same bias (which is often the case), (3.1) implies that variance estimation with OL subsamples is asymptotically more efficient than the NOL subsample alternative owing to a smaller asymptotic MSE.

Unlike $K_0$, $K_1$ does depend on the volume $|R_0|$, which in turn is constrained by the $R_0$-template's geometry. Through $|R_0|$ in (3.1), $K_1$ is ultimately bounded by the amount of space that an object of $R_0$'s *shape* can possibly occupy within $(-1/2, 1, 2]^d$ [i.e., by how much volume can be filled by a given geometrical body (e.g., circle) compared to a cube]. The constants $K_1$ in Table 1 are computed with templates of prescribed shape and largest possible volume in $(-1/2, 1/2]^d$. These values most accurately reflect the influence of $R_0$'s (or $R_n$'s) geometry on the large-sample relative performance of $\hat{\tau}_{n,\mathrm{OL}}^2$ and $\hat{\tau}_{n,\mathrm{NOL}}^2$ in terms of variance in (3.1) and also efficiency (see Section 5).

To conclude this section, we remark that both subsample-based variance estimators can be shown to be MSE-consistent under Theorem 3.1 conditions, allowing for more general spatial sampling regions, in both shape and dimension, than previously considered. Inference on the parameter $\theta$ can be made through the limiting standard normal distribution of $\sqrt{N_n}(\hat{\theta}_n - \theta)/\hat{\tau}_n$ for $\hat{\tau}_n = \hat{\tau}_{n,\mathrm{OL}}$ or $\hat{\tau}_{n,\mathrm{NOL}}$.

**4. Bias expansions.** We now try to capture and precisely describe the leading order terms in the asymptotic bias of each subsample-based variance estimator, similar to the variance determinations from the previous section. We first establish and note the order of the dominant component in the bias expansions of $\hat{\tau}_{n,\mathrm{OL}}^2$ and $\hat{\tau}_{n,\mathrm{NOL}}^2$, which is the subject of the following lemma.

LEMMA 4.1. *With Assumptions* A.1–A.5, *suppose that Conditions $D_2$ and $M_{2+a}$ hold for $d \geq 2$ or that $D_3$ and $M_{3+a}$ hold for $d = 1$ (where a is as specified by the respective Condition $D_r$). Then the subsample estimators of*

$\tau_n^2 = N_n \operatorname{Var}(\hat{\theta}_n)$ *have expectations*

$$\operatorname{E}(\hat{\tau}_{n,\mathrm{OL}}^2) = \tau_n^2 + O(1/{}_s\lambda_1^{(n)}) \quad and \quad \operatorname{E}(\hat{\tau}_{n,\mathrm{NOL}}^2) = \tau_n^2 + O(1/{}_s\lambda_1^{(n)}).$$

The lemma shows that, under the smooth function model, the asymptotic bias of each estimator is $O(1/{}_s\lambda_1^{(n)})$ for all dimensions of sampling. Politis and Romano (1993a) and Sherman (1996) showed this same *size* for the bias of $\hat{\tau}_{n,\mathrm{OL}}^2$ with sampling regions based on rectangles $R_0 = (-1/2, 1/2]^d$ or simple closed curves in $\mathbb{R}^2$, respectively. Lemma 4.1 extends these results to a broader class of sampling regions. However, we would like to precisely identify the $O(1/{}_s\lambda_1^{(n)})$ bias component for $\hat{\tau}_{n,\mathrm{OL}}^2$ or $\hat{\tau}_{n,\mathrm{NOL}}^2$ to obtain optimal subsample scaling that accounts for the geometry of $R_n$.

To achieve some measure of success in determining the exact bias of the subsampling estimators, we reformulate the subsampling design slightly so that ${}_s\lambda_n \equiv {}_s\lambda_1^{(n)} = \cdots = {}_s\lambda_d^{(n)}$. That is, a common scaling factor in all directions is now used to define the subsampling regions, as in Sherman and Carlstein (1994) and Sherman (1996). This constraint will allow us to deal with the counting issues at the heart of the bias expansion.

Adopting a common scaling factor ${}_s\lambda_n$ for the subsamples also is sensible for a few other reasons at this stage:

1. "Unconstrained" optimum values of ${}_s\Delta_n$ cannot always be found by minimizing the asymptotic MSE of $\hat{\tau}_{n,\mathrm{OL}}^2$ or $\hat{\tau}_{n,\mathrm{NOL}}^2$, even for variance estimation of some desirable statistics on geometrically "simple" sampling and subsampling regions. Consider estimating the variance of a real-valued sample mean over a rectangular sampling region in $\mathbb{R}^d$ based on $R_0 = (-1/2, 1/2]^d$, with observations on $\mathbf{Z}^d = \mathbb{Z}^d$. If Assumptions A.1–A.5 and Condition $M_1$ hold, the leading term in the bias expansion can be shown to be

$$\text{Bias of } \hat{\tau}_{n,\mathrm{OL}}^2 = \left(-\sum_{i=1}^{d} \frac{L_i}{{}_s\lambda_i^{(n)}}\right)(1 + o(1));$$

$$L_i = \sum_{\substack{\mathbf{k} \in \mathbb{Z}^d \\ \mathbf{k} = (k_1, \ldots, k_d)'}} |k_i| \operatorname{Cov}(Z(\mathbf{0}), Z(\mathbf{k})).$$

In using the parenthetical sum above to expand the MSE of $\hat{\tau}_{n,\mathrm{OL}}^2$, one finds that the resulting MSE cannot be minimized over the permissible, positive range of ${}_s\Delta_n$ if the signs of the $L_i$ values are unequal. That is, for $d > 1$, the subsample estimator MSE cannot always be globally minimized to obtain optimal subsample factors ${}_s\Delta_n$ by considering just the leading order bias terms. An effort to determine and incorporate (into the asymptotic MSE) second- or third-order bias components quickly becomes intractable, even with rectangular regions.

2. The diagonal components of $_s\Delta_n$ are asymptotically scalar multiples of each other by Assumption A.1. If so desired, a template choice for $R_0$ could be used to scale the expansion of the subsampling regions in each direction.

In the continuing discussion, we assume

$$_sR_n = {_s\lambda_n} R_0. \tag{4.1}$$

We frame the components necessary for determining the biases of the spatial subsample variance estimators in the next theorem. Let

$$C_n(\mathbf{k}) \equiv |\mathbf{Z}^d \cap {_sR_n} \cap (\mathbf{k} + {_sR_n})|$$

denote the number of pairs of observations in the subsampling region $_sR_n$ separated by a translate $\mathbf{k} \in \mathbb{Z}^d$.

THEOREM 4.1. *Suppose that $d \geq 2$, $_sR_n = {_s\lambda_n} R_0$ and Assumptions A.1–A.5, Conditions $D_3$ and $M_{3+a}$ hold with $a$ as specified under Condition $D_3$. If, in addition, $_s\lambda_n \in \mathbb{Z}_+$ for NOL subsamples and*

$$\lim_{n \to \infty} \frac{_sN_n - C_n(\mathbf{k})}{({_s\lambda_n})^{d-1}} = C(\mathbf{k}) \tag{4.2}$$

*exists for all $\mathbf{k} \in \mathbb{Z}^d$, then*

$$\mathrm{E}(\hat{\tau}_n^2) - \tau_n^2 = \frac{-1}{_s\lambda_n |R_0|} \left( \sum_{\mathbf{k} \in \mathbb{Z}^d} C(\mathbf{k}) \sigma(\mathbf{k}) \right) (1 + o(1)),$$

*where $\sigma(\mathbf{k}) = \mathrm{Cov}(\nabla' Z(\mathbf{t}), \nabla' Z(\mathbf{t} + \mathbf{k}))$ and where $\hat{\tau}_n^2$ is either $\hat{\tau}_{n,\mathrm{OL}}^2$ or $\hat{\tau}_{n,\mathrm{NOL}}^2$.*

Note that the numerator on the left-hand side of (4.2) is the number of $\mathbf{Z}^d$ grid points that lie in the subregion $_sR_n$, but not in the translate $\mathbf{k} + {_sR_n}$. Hence, computing the bias above actually requires counting the number of lattice points inside intersections like $_sR_n \cap \mathbf{k} + {_sR_n}$, which is difficult in general. To handle the problem, one may attempt to estimate the count $C_n(\mathbf{k})$ with the corresponding Lebesgue volume, $|_sR_n \cap \mathbf{k} + {_sR_n}|$, and then quantify the resulting approximation error. The determination of volumes or areas may not be easy either but hopefully more manageable. For example, if $R_0$ is a circle, the area of $_s\lambda_n R_0$ can be readily computed, but the number of $\mathbb{Z}^2$ integers inside $_s\lambda_n R_0$ is not so simple and was in fact a famous consideration of Gauss [cf. Krätzel (1988), page 141].

We first note that the boundary condition on $R_0$ provides a general (trivial) bound on the discrepancy between the count $C_n(\mathbf{k})$ and the volume $|_sR_n \cap \mathbf{k} + {_sR_n}|$: $O(_s\lambda_n{}^{d-1})$. However, the size of the numerator in (4.2) is also $O(_s\lambda_n{}^{d-1})$, corresponding to the order of $\mathbf{Z}^d$ lattice points "near" the boundary of ${}_sR_n$. Consequently, a standard $O({}_s\lambda_n{}^{d-1})$ bound on the volume-count approximation error is too large to immediately justify the exchange of volumes $|{}_sR_n|$, $|{}_sR_n \cap \mathbf{k} + {}_sR_n|$ for counts ${}_sN_n$, $C_n(\mathbf{k})$ in (4.2).

Bounds on the difference between lattice point counts and volumes have received much attention in analytic number theory, which we briefly mention. Research has classically focused on sets outlined by "smooth" simple closed curves in the plane $\mathbb{R}^2$ and on one question in particular [Huxley (1996)]: When a curve with interior area $A$ is "blown up" by a factor $b$, how large is the difference between the number of $\mathbb{Z}^2$ integer points inside the new curve and the area $b^2A$? For convex sets with a *smoothly winding border*, van der Corput's (1920) answer to the posed question above is $O(b^{46/69+\varepsilon})$, while the best answer is $O(b^{46/73+\varepsilon})$ for curves with sufficiently *differentiable radius of curvature* [Huxley (1993, 1996)]. These types of bounds, however, are invalid for many convex polygonal templates $R_0$ in $\mathbb{R}^2$ such as triangles, trapezoids, and so on, where often the difference between number of $\mathbb{Z}^2$ integer points in ${}_sR_n = {}_s\lambda_n R_0$ and its area is of exact order $O({}_s\lambda_n)$ (set also by the boundary condition on $R_0$ or the perimeter length of ${}_sR_n$). The problem above, as considered by number theorists, does not directly address counts for intersections between an expanding region and its vector translates, for example, ${}_sR_n \cap \mathbf{k} + {}_sR_n$.

To eventually compute closed-form bias expansions for $\hat{\tau}_{n,\text{OL}}^2$, we use approximation techniques for *subtracted* lattice point counts. For each $\mathbf{k} \in \mathbb{Z}^d$, we:

1. Replace the numerator of (4.2) with the difference of corresponding Lebesgue volumes.
2. Show the following error term is of sufficiently small order $o({}_s\lambda_n^{d-1})$:

$$({}_sN_n - C_n(\mathbf{k})) - ({}_s\lambda_n^d|R_0| - |{}_sR_n \cap \mathbf{k} + {}_sR_n|)$$
$$= ({}_sN_n - {}_s\lambda_n^d|R_0|) - (C_n(\mathbf{k}) - |{}_sR_n \cap \mathbf{k} + {}_sR_n|).$$

We do approximate the number of lattice points in ${}_sR_n$ and ${}_sR_n \cap \mathbf{k} + {}_sR_n$ by set volumes, though the Lebesgue volume may not adequately capture the lattice point count in either set. However, the *difference* between approximation errors ${}_sN_n - {}_s\lambda_n^d|R_0|$ and $C_n(\mathbf{k}) - |{}_sR_n \cap \mathbf{k} + {}_sR_n|$ can be shown to be asymptotically small enough, for some templates $R_0$, to justify replacing counts with volumes in (4.2) (see Lemma 10.4). That is, these two volume count estimation errors can cancel to a sufficient extent when subtracted. The above approach becomes slightly more complicated for NOL subsamples, $\tilde{R}_{\mathbf{i},n} = {}_s\Delta_n(\mathbf{i} + R_0)$, which may vary in number of sampling sites ${}_sN_{\mathbf{i},n}$. In this case, errors incurred by approximating counts $|\mathbf{Z}^d \cap \tilde{R}_{\mathbf{i},n} \cap \mathbf{k} + \tilde{R}_{\mathbf{i},n}|$ with volumes $|\tilde{R}_{\mathbf{i},n} \cap \mathbf{k} + \tilde{R}_{\mathbf{i},n}|$ are shown to be asymptotically negligible, uniformly in $\mathbf{i} \in J_{\text{NOL}}$.

In the following theorem, we use this technique to give bias expansions for a large class of sampling regions in $\mathbb{R}^d$, $d \leq 3$, which are "nearly" convex. The sampling region $R_n$ may differ from a convex set possibly only at its boundary, but sampling sites on the border may be arbitrarily included or excluded from $R_n$.

Some notation is additionally required. For $\alpha = (\alpha_1, \ldots, \alpha_p)' \in (\mathbb{Z}_+)^p$, $\mathbf{x} \in \mathbb{R}^p$, write $\mathbf{x}^\alpha = \prod_{i=1}^p x_i^{\alpha_i}$, $\alpha! = \prod_{i=1}^p (\alpha_i!)$, and $c_\alpha = D^\alpha H(\mu)/\alpha!$. Let $Z_\infty$ denote a random vector with a normal $\mathcal{N}(0, \Sigma_\infty)$ distribution on $\mathbb{R}^p$, where $\Sigma_\infty$ is the limiting covariance matrix of the scaled sample mean $\sqrt{N_n}(\bar{Z}_{N_n} - \mu)$ from (2.1). Let $B^\circ$, $\overline{B}$ denote the interior and closure of $B \subset \mathbb{R}^d$, respectively.

THEOREM 4.2. *Suppose ${}_sR_n = {}_s\lambda_n R_0$ and there exists a convex set $B$ such that $B^\circ \subset R_0 \subset \overline{B}$. With Assumptions A.2–A.5, assume Conditions $D_{5-d}$ and $M_{5-d+a}$ hold for $d \in \{1, 2, 3\}$ (where $a$ is as specified by the respective Condition $D_r$). Then*

$$C(\mathbf{k}) = V(\mathbf{k}) \equiv \lim_{n \to \infty} \frac{|{}_sR_n| - |{}_sR_n \cap (\mathbf{k} + {}_sR_n)|}{({}_s\lambda_n)^{d-1}}, \qquad \mathbf{k} \in \mathbb{Z}^d,$$

*whenever $V(\mathbf{k})$ exists and the biases $\mathrm{E}(\hat{\tau}^2_{n,\mathrm{OL}}) - \tau_n^2$, $\mathrm{E}(\hat{\tau}^2_{n,\mathrm{NOL}}) - \tau_n^2$ are equal to, for $d = 1$,*

$$\frac{-1}{{}_s\lambda_n |R_0|} \left( \sum_{\mathbf{k} \in \mathbb{Z}} |\mathbf{k}| \sigma(\mathbf{k}) + C_\infty \right)(1 + o(1));$$

*for $d = 2$ or $3$,*

$$\left( -\sum_{\mathbf{k} \in \mathbb{Z}^d} \frac{|{}_sR_n| - |{}_sR_n \cap (\mathbf{k} + {}_sR_n)|}{|{}_sR_n|} \sigma(\mathbf{k}) \right)(1 + o(1))$$

*or*

$$\frac{-1}{{}_s\lambda_n |R_0|} \left( \sum_{\mathbf{k} \in \mathbb{Z}^d} V(\mathbf{k}) \sigma(\mathbf{k}) \right)(1 + o(1)),$$

*provided each $V(\mathbf{k})$ exists, where $\sigma(\mathbf{k}) = \mathrm{Cov}(\nabla' Z(\mathbf{t}), \nabla' Z(\mathbf{t} + \mathbf{k}))$ and*

$$C_\infty = \mathrm{Var}\left( \sum_{\|\alpha\|_1 = 2} \frac{c_\alpha}{\alpha!} Z_\infty^\alpha \right) + 2 \sum_{\substack{\|\alpha\|_1 = 1 \\ \|\beta\|_1 = 3}} \frac{c_\alpha c_\beta}{\beta!} \mathrm{E}(Z_\infty^\alpha Z_\infty^\beta)$$

$$+ 2 \sum_{\mathbf{k}_1, \mathbf{k}_2 \in \mathbb{Z}} \sum_{\substack{\|\alpha\|_1 = 1 \\ \|\beta\|_1 = 1, \|\gamma\|_1 = 1}} \frac{c_\alpha c_{(\beta + \gamma)}}{(\beta + \gamma)!}$$

$$\times \mathrm{E}([Z(\mathbf{t}) - \mu]^\alpha [Z(\mathbf{t} + \mathbf{k}_1) - \mu]^\beta [Z(\mathbf{t} + \mathbf{k}_2) - \mu]^\gamma).$$

REMARK 4.1. If Condition $D_m$ holds with $\mathcal{C}=0$ for some $m \in \{2,3,4\}$, then Condition $M_{m-1}$ is sufficient in Theorem 4.2.

REMARK 4.2. For each $\mathbf{k} \in \mathbb{Z}^d$, the numerator in $V(\mathbf{k})$ is $O({_s\lambda_n}^{d-1})$ by the $R_0$-boundary condition which holds for convex templates. We may then expand the bias of the estimators through the limiting, scaled volume differences $V(\mathbf{k})$. For $d = 1$, with samples and subsamples based on intervals, it can be easily seen that $V(\mathbf{k}) = |\mathbf{k}|$, which appears in Theorem 4.2.

The function $H(\cdot)$ needs to be increasingly "smoother" to determine the bias component of $\hat{\tau}^2_{n,\mathrm{OL}}$ or $\hat{\tau}^2_{n,\mathrm{NOL}}$ in lower-dimensional spaces $d = 1$ or $2$. For a real-valued time series sample mean $\hat{\theta}_n = \bar{Z}_n$, the well-known bias of the subsample variance estimators follows from Theorem 4.2 under our sampling framework $R_0 = (-1/2, 1/2]$, $\mathbf{Z} = \mathbb{Z}$ as

$$(4.3) \qquad \frac{-1}{{_s\lambda_n}} \left( \sum_{\mathbf{k} \in \mathbb{Z}} |\mathbf{k}| \operatorname{Cov}(\nabla' Z(\mathbf{0}), \nabla' Z(\mathbf{k})) \right)$$

with $\nabla = 1$. In general though, terms in the Taylor expansion of $\hat{\theta}_{\mathbf{i},n}$ (around $\mu$) up to fourth order can contribute to the bias of $\hat{\tau}^2_{n,\mathrm{OL}}$ and $\hat{\tau}^2_{n,\mathrm{NOL}}$ when $d = 1$. In contrast, the asymptotic bias of the time series MBB variance estimator with "smooth" model statistics is very different from its subsample-based counterpart. The MBB variance estimator's bias is given by (4.3), determined only by the linear component from the Taylor expansion of $\hat{\theta}_{\mathbf{i},n}$ [cf. Lahiri (1996)].

**5. Asymptotically optimal subsample sizes.** In the following, we consider "size" selection for the subsampling regions to maximize the large-sample accuracy of the subsample variance estimators. For reasons discussed in Section 4, we examine a theoretically optimal scaling choice ${_s\lambda_n}$ for subregions in (4.1).

5.1. *Theoretical optimal subsample sizes.* Generally speaking, there is a trade-off in the effect of subsample size on the bias and variance of $\hat{\tau}^2_{n,\mathrm{OL}}$ or $\hat{\tau}^2_{n,\mathrm{NOL}}$. Increasing ${_s\lambda_n}$ reduces the bias but increases the variance of the estimators. The best value of ${_s\lambda_n}$ optimizes the overall performance of a subsample variance estimator by balancing the contributions from both the estimator's variance and bias. An optimal ${_s\lambda_n}$ choice can be found by minimizing the asymptotic order of a variance estimator's MSE under a given OL or NOL sampling scheme.

Theorem 4.1 implies that the bias of the estimators $\hat{\tau}^2_{n,\mathrm{OL}}$ and $\hat{\tau}^2_{n,\mathrm{NOL}}$ is of exact order $O(1/{_s\lambda_n})$. For a broad class of sampling regions $R_n$, the leading order bias component can be determined explicitly with Theorem 4.2. We bring these variance and bias expansions together to obtain an optimal subsample scaling factor $_s\lambda_n^{\text{opt}}$.

THEOREM 5.1. *Let $_sR_n = {_s\lambda_n}R_0$. With Assumptions A.2–A.5, assume Conditions $D_2$ and $M_{5+2a}$ hold if $d \geq 2$ or Conditions $D_3$ and $M_{7+2a}$ hold if $d = 1$ (where $a$ is as specified by the respective Condition $D_r$). If*

$$B_0|R_0| \equiv \sum_{\mathbf{k} \in \mathbb{Z}^d} C(\mathbf{k})\sigma(\mathbf{k}) + I_{\{d=1\}}C_\infty \neq 0,$$

*then*

$$_s\lambda_{n,\text{OL}}^{\text{opt}} = \left(\frac{\det(\Delta_n)(B_0)^2}{dK_0\tau^4}\right)^{1/(d+2)}(1+o(1))$$

*and*

$$_s\lambda_{n,\text{NOL}}^{\text{opt}} = \left(\frac{\det(\Delta_n)|R_0|(B_0)^2}{d\tau^4}\right)^{1/(d+2)}(1+o(1)).$$

REMARK 5.1. If Condition $D_m$ holds with $\mathcal{C} = 0$ for some $m \in \{2,3\}$, then Condition $M_{2m-1}$ is sufficient.

REMARK 5.2. Theorem 5.1 suggests that optimally scaled OL subsamples should be *larger* than the NOL ones by a scalar: $(K_1)^{-1/(d+2)} > 1$ where $K_1 = K_0|R_0|$ is the limiting ratio of variances from (3.1).

It is well known in the time series case that the OL subsampling scheme produces an asymptotically more efficient variance estimator than its NOL counterpart. We can now quantify the relative efficiency of the two subsampling procedures in $d$-dimensional sampling space. With each variance estimator respectively optimized using (4.1), $\hat{\tau}_{n,\text{OL}}^2$ is more efficient than $\hat{\tau}_{n,\text{NOL}}^2$ and the asymptotic relative efficiency $(ARE_d)$ of $\hat{\tau}_{n,\text{NOL}}^2$ to $\hat{\tau}_{n,\text{OL}}^2$ depends solely on the geometry of $R_0$,

$$ARE_d = \lim_{n \to \infty} \frac{\text{E}(\hat{\tau}_{n,\text{OL}}^2 - \tau_n^2)^2}{\text{E}(\hat{\tau}_{n,\text{NOL}}^2 - \tau_n^2)^2} = (K_1)^{2/(d+2)} < 1.$$

Possolo (1991), Politis and Romano (1993a, 1994), Hall and Jing (1996) and Garcia-Soidan and Hall (1997) have examined subsampling with rectangular regions based essentially on $R_0 = (-1/2, 1/2]^d$. Using the geometrical characteristic $K_1 = (\frac{2}{3})^d$ for rectangles, we can now examine the effect of the sampling dimension on the $ARE_d$ of $\hat{\tau}_{n,\text{NOL}}^2$ to $\hat{\tau}_{n,\text{OL}}^2$ for these sampling regions. Although the $ARE_d$ decreases as the dimension $d$ increases, we find the relative improvement of $\hat{\tau}_{n,\text{OL}}^2$ over $\hat{\tau}_{n,\text{NOL}}^2$ is ultimately limited and the $ARE_d$ has a lower bound of $4/9$ for all $\mathbb{R}^d$-rectangular regions.

5.2. *Theoretical optimal subsample shapes.* We conclude this section by addressing a question raised by a referee on subsample *shape* selection. Although not widely considered in the literature, subsample variance estimators are also possible by using subsamples of a freely chosen shape, rather than scaled-down copies of $R_n$. Nordman and Lahiri (2003) discuss comparing variance estimators, based on differently shaped subsamples, through their asymptotic relative efficiency. This involves finding MSE expansions for estimators with OL, NOL subsamples of an arbitrary shape with optimal scaling (e.g., modified versions of Theorems 3.1, 4.1 and 5.1). However, because both the subsample geometry and the r.f. covariances influence a subsample estimator's bias (see Section 6), a direct comparison of asymptotic MSEs to choose an optimal subsample shape can become complicated, especially for OL subsamples.

For illustration, consider selecting between circular and rectangular subsamples for sample mean $\hat{\theta}_n = \bar{Z}_{N_n} \in \mathbb{R}$ variance estimation on a rectangular region $R_n \subset \mathbb{R}^2$ under a Gaussian isotropic covariogram,

$$\sigma(\mathbf{k}) = \exp(-\beta \|\mathbf{k}\|^2), \qquad \mathbf{k} \in \mathbb{Z}^2.$$

The value of $\beta$ heavily affects the large sample performances of circles and rectangles (i.e., scaled-down copies of $R_n$) as subsamples and makes the choice of subsample shape difficult. For example, the asymptotic efficiency of circular to rectangular OL (NOL) subsamples is 0.9259 (1.0274) for $\beta = 0.2$ and 1.0758 (1.1937) for $\beta = 2$. We conducted a small simulation study of the finite sample efficiencies of these subsample shapes on several rectangular $R_n$ to compare with the asymptotic values. The results in Table 4 indicate that the asymptotic advantages of a subsample shape may also not be readily apparent in finite samples due to edge effects. See Nordman and Lahiri (2003) for further details and examples on the effect of subsample shape for variance estimation.

**6. Examples.** We now provide some examples of the important quantities $K_0$, $K_1$, $B_0$ associated with optimal scaling $_s\lambda_n^{\text{opt}}$ with some common sampling region templates, determined from Theorems 3.1 and 4.2. For subsamples from (4.1), the theoretically best $_s\lambda_n^{\text{opt}}$ can also be formulated in terms of $|R_n| = \det(\Delta_n)|R_0|$ (sampling region volume), $K_1$ and $B_0$.

6.1. *Examples in* $\mathbb{R}^2$.

EXAMPLE 1. Rectangular regions in $\mathbb{R}^2$ (potentially rotated): if

$$R_0 = \{((l_1 \cos\theta, l_2 \sin\theta)\mathbf{x}, (-l_1 \sin\theta, l_2 \cos\theta)\mathbf{x})' : \mathbf{x} \in (-1/2, 1/2]^2\}$$

TABLE 4
*Minimal normalized MSE* $\mathrm{E}(\hat{\tau}_n^2/\tau_n^2 - 1)^2$ *for OL/NOL subsample estimators* $\hat{\tau}_n^2$ *of sample mean variance* $\tau_n^2 = N_n \operatorname{Var}(\bar{Z}_{N_n})$ *on* $R_n \cap \mathbb{Z}^2$, *with* $\sigma(\mathbf{k}) = \exp(-\beta\|\mathbf{k}\|^2)$, $\mathbf{k} \in \mathbb{Z}^2$ *(based on* 1000 *simulations). Rectangular (rec.) and circular (cir.) subsamples* $_s\lambda_n^{\mathrm{opt}} R_0^*$ *are based on* $R_0^* = (-1/2, 1/2]^2$, $\{\mathbf{x} \in \mathbb{R}^2 : \|\mathbf{x}\| \leq 1/2\}$ *using optimal scaling* $_s\lambda_n^{\mathrm{opt}}$ *(an integer listed beside each MSE). Estimated relative efficiencies (RE) of cir. versus rec. subsamples are also listed*

| $R_n$ | rec. subsamples | | cir. subsamples | | cir./rec. $RE$ | |
|---|---|---|---|---|---|---|
| | OL | NOL | OL | NOL | OL | NOL |
| | | | $\beta = 0.2$ | | | |
| $(-5, 5]^2$ | 0.4295 (4) | 0.4261 (5) | 0.4519 (2) | 0.4286 (5) | 1.0521 | 1.0060 |
| $(-10, 10]^2$ | 0.2329 (5) | 0.2183 (5) | 0.2418 (5) | 0.2328 (5) | 1.0384 | 1.0661 |
| $(-30, 30]^2$ | 0.0806 (10) | 0.0842 (10) | 0.0835 (9) | 0.0944 (9) | 1.0355 | 1.1260 |
| $(-50, 50]^2$ | 0.0482 (14) | 0.0562 (11) | 0.0462 (15) | 0.0601 (11) | 0.9585 | 1.0698 |
| | | | $\beta = 2$ | | | |
| $(-5, 5]^2$ | 0.0841 (2) | 0.0978 (2) | 0.1170 (2) | 0.1426 (1) | 1.3890 | 1.4570 |
| $(-10, 10]^2$ | 0.0436 (3) | 0.0515 (2) | 0.0436 (3) | 0.0641 (3) | 1.0000 | 1.2432 |
| $(-30, 30]^2$ | 0.0128 (5) | 0.0162 (4) | 0.0138 (5) | 0.0199 (5) | 1.0771 | 1.2260 |
| $(-50, 50]^2$ | 0.0082 (6) | 0.0111 (5) | 0.0092 (7) | 0.0129 (5) | 1.1139 | 1.1594 |

for $\theta \in [0, \pi]$, $0 < l_1, l_2$, then

$$K_0 = \frac{4}{9}, \qquad B_0 = \sum_{\substack{\mathbf{k} \in \mathbb{Z}^2 \\ \mathbf{k} = (k_1, k_2)'}} \left( \frac{|k_1 \cos\theta - k_2 \sin\theta|}{l_1} + \frac{|k_1 \sin\theta + k_2 \cos\theta|}{l_2} \right) \sigma(\mathbf{k}).$$

The characteristics $K_1$, $B_0$ for determining optimal subsamples based on two rectangular templates, including a diamond-shaped region (i.e., $\theta = \pi/4$, $l_1 = l_2 = 1/\sqrt{2}$), are further described in Table 5.

EXAMPLE 2. If $R_0$ is a circle of radius $r \leq 1/2$ centered at the origin, then $K_0$ appears in Table 3 and $B_0 = 2/(r\pi) \sum_{\mathbf{k} \in \mathbb{Z}^2} \|\mathbf{k}\| \sigma(\mathbf{k})$.

EXAMPLE 3. For *any* triangle, $K_0 = 2/5$. Two examples are provided in Tables 2 and 5.

EXAMPLE 4. If $R_0$ is a regular hexagon, centered at the origin and with side length $l \leq 1/2$, then

$$K_0 = \frac{37}{81}, \qquad B_0 = \frac{2\sqrt{3}}{l} \sum_{\mathbf{k} \in \mathbb{Z}^2} (|k_2| + \max\{\sqrt{3}|k_1|, |k_2|\}) \sigma(\mathbf{k}).$$

TABLE 5
*Examples of several shapes of $R_0 \subset \mathbb{R}^2$ and associated $K_1$, $B_0$ for ${}_s\lambda_n^{\mathrm{opt}}$*

| $R_0$ | $K_1$ | $B_0$ |
|---|---|---|
| $(-1/2, 1/2]^2$ | $4/9$ | $\sum_{\mathbf{k} \in \mathbb{Z}^2} \|\mathbf{k}\|_1 \sigma(\mathbf{k})$ |
| Circle of radius $1/2$ at origin | $\pi/4 - 4/(3\pi)$ | $4/\pi \sum_{\mathbf{k} \in \mathbb{Z}^2} \|\mathbf{k}\| \sigma(\mathbf{k})$ |
| Diamond in Figure 1(i) | $2/9$ | $2 \sum_{\mathbf{k} \in \mathbb{Z}^2} \|\mathbf{k}\|_\infty \sigma(\mathbf{k})$ |
| Right triangle in Figure 1(ii) | $1/5$ | Table 2 |
| Triangle in Figure 1(iii) | $1/5$ | $\sum_{\mathbf{k} \in \mathbb{Z}^2} (|k_2| + \max\{2|k_1|, |k_2|\}) \sigma(\mathbf{k})$ |
| Parallelogram in Figure 1(iv) | $2/9 + (\sqrt{5}-1)/375$ | $4/\sqrt{5} \sum_{\mathbf{k} \in \mathbb{Z}^2} (|k_1 - 2k_2|/5 + |k_2|) \sigma(\mathbf{k})$ |

EXAMPLE 5. For any parallelogram in $\mathbb{R}^2$ with interior angle $\gamma$ and adjacent sides of ratio $b \geq 1$, $K_0 = 4/9 + 2/15 \cdot b^{-2} |\cos\gamma| (1 - |\cos\gamma|)$. In particular, if a parallelogram $R_0$ is formed by two vectors $(0, l_1)'$, $(l_2 \cos\gamma, l_2 \sin\gamma)'$ extended from a point $\mathbf{x} \in (-1/2, 1/2]^2$, then

$$B_0 = \frac{1}{|\sin\theta|} \sum_{\mathbf{k} \in \mathbb{Z}^2} \left( \frac{|k_1 \cdot |\cos\theta| - k_2 \cdot |\sin\theta||}{\max\{l_1, l_2\}} + \frac{|k_2|}{\min\{l_1, l_2\}} \right) \sigma(\mathbf{k}),$$

$$\gamma \in (0, \pi),\ l_1, l_2 > 0.$$

For further bias term $B_0$ calculation tools with more general (nonconvex) sampling regions and templates $R_0$ (represented as the union of two approximately convex sets), see Nordman (2002).

6.2. *Examples in $\mathbb{R}^d$, $d \geq 3$.*

EXAMPLE 6. For any sphere, $K_0$ is given in Table 3. The properties $B_0$, $K_1$ of the sphere described in Tables 1 and 2 correspond to the template sphere $R_0$ of radius $1/2$ with maximal volume in $(-1/2, 1/2]^3$.

EXAMPLE 7. The $K_0$ value for any $\mathbb{R}^3$ cylinder appears in Table 3. If $R_0$ is a cylinder with circular base (parallel to the x–y plane) of radius $r$

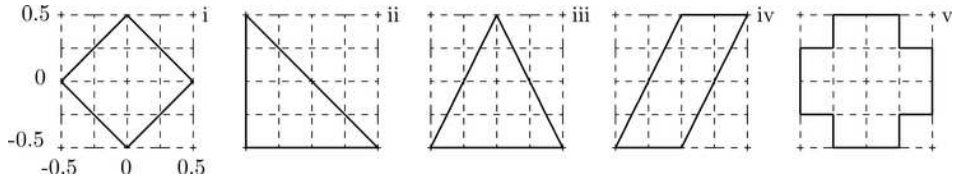

FIG. 1. *Examples of templates $R_0 \subset (-1/2, 1/2]^2$ are outlined by solid lines. Cross-shaped sampling regions $R_n$ described in Table 2 are based on $R_0$ in* (v).

and height $h$, then

$$B_0 = \sum_{\substack{\mathbf{k} \in \mathbb{Z}^3 \\ \mathbf{k}=(k_1,k_2,k_3)'}} \left( \frac{|k_3|}{h} + \frac{2\sqrt{k_1^2 + k_2^2}}{\pi r} \right) \sigma(\mathbf{k}).$$

The results of Theorem 4.2 for determining the bias $B_0$ also seem plausible for convex sampling regions in $\mathbb{R}^d$, $d \geq 4$, but require further study of lattice point counting techniques in higher dimensions. However, bias expansions of the OL and NOL subsample variance estimators are relatively straightforward for an important class of rectangular sampling regions based on the prototype $R_0 = (-1/2, 1/2]^d$, which can then be used in optimal subsample scaling. These hypercubes have "faces" parallel to the coordinate axes, which simplifies the task of counting sampling sites, or lattice points, within such regions. We give precise bias expansions in the following theorem, while allowing for potentially missing sampling sites at the border of the sampling region $R_n$.

THEOREM 6.1. *Let $(-1/2, 1/2]^d \subset \Lambda_\ell^{-1} R_0 \subset [-1/2, 1/2]^d$, $d \geq 3$, for a $d \times d$ diagonal matrix $\Lambda_\ell$ with entries $0 < \ell_i \leq 1$, $i = 1, \ldots, d$. Suppose $_sR_n = {}_s\lambda_n R_0$ and Assumptions A.2–A.5, Conditions $D_2$ and $M_{2+a}$ hold with a as specified under Condition $D_2$. Then the biases $\mathrm{E}(\hat{\tau}_{n,\mathrm{OL}}^2) - \tau_n^2$, $\mathrm{E}(\hat{\tau}_{n,\mathrm{NOL}}^2) - \tau_n^2$ are equal to $-{}_s\lambda_n^{-1} B_0(1 + o(1))$ where*

$$B_0 = \sum_{\mathbf{k} \in \mathbb{Z}^d} \left( \sum_{i=1}^d \frac{|k_i|}{\ell_i} \right) \sigma(\mathbf{k}), \qquad \sigma(\mathbf{k}) = \mathrm{Cov}(\nabla' Z(\mathbf{t}), \nabla' Z(\mathbf{t} + \mathbf{k})).$$

EXAMPLE 8. For rectangular sampling regions $R_n = \Delta_n(-1/2, 1/2]^d$, optimal subsamples (4.1) may be chosen with

$$_s\lambda_{n,\mathrm{NOL}}^{\mathrm{opt}} = \left( \frac{|R_n|}{d\tau^4} \left( \sum_{\mathbf{k} \in \mathbb{Z}^d} \|\mathbf{k}\|_1 \sigma(\mathbf{k}) \right)^2 \right)^{1/(d+2)} (1 + o(1))$$

or

$$_s\lambda_{n,\mathrm{OL}}^{\mathrm{opt}} = {}_s\lambda_{n,\mathrm{NOL}}^{\mathrm{opt}} (\tfrac{3}{2})^{d/(d+2)},$$

using the template $R_0 = (-1/2, 1/2]^d$.

**7. Empirical subsample size determination.** This section considers data-based estimation of the theoretical optimal scaling factor $_s\lambda_n^{\mathrm{opt}}$ for subsamples as in (4.1). We describe two estimation techniques for this. One approach involves using "plug-in" estimates and the second involves minimizing an estimated MSE criterion function. In Section 8 we evaluate both estimation methods for $_s\lambda_n^{\text{opt}}$ through a simulation study. Inference on "best" subsample scaling closely resembles the problem of empirically gauging the theoretically optimal block length with the MBB variance estimator. With time series, estimation rules of optimal MBB block size have been developed using both plug-in and empirical MSE methods [cf. Bühlmann and Künsch (1999) and Hall, Horowitz and Jing (1995)].

Hall and Jing (1996) give a method for estimating optimal subsample scaling through minimization of an estimated MSE function in the time series case. Considering OL subsamples first, we adapt this approach (hereafter the HJ method) for spatial subsampling as follows. We determine the template $R_0$ as the largest set of the form $\Delta_n^{-1} R_n$ within $(-1/2, 1/2]^d$. Let $J_{\text{OL}}(\lambda_m)$ denote a collection of OL subsamples using a scaling factor $_s\lambda_n \equiv \lambda_m > 0$ in (4.1). Here $\lambda_m$ is a "smoothing parameter." We treat each subsample in $J_{\text{OL}}(\lambda_m)$ as a scale $\lambda_m R_0$ sampling region on which an OL subsample variance estimator, with subsample scaling $_s\lambda_m < \lambda_m$, can be computed. Denote the resulting variance estimates as $\hat{\tau}_{i,m,\text{OL}}^2$, $i = 1, \ldots, |J_{\text{OL}}(\lambda_m)|$. Write $\hat{\tau}_{n,\text{OL}}^2 \equiv \hat{\tau}_{n,\text{OL}}^2(\lambda_m)$ as the variance estimator computed on the region $R_n$ with subsample scaling $\lambda_m$. An estimate of the MSE when using subsamples of size $_s\lambda_m R_0$ on regions of size $\lambda_m R_0$ is the average of the squared differences $(\hat{\tau}_{i,m,\text{OL}}^2 - \hat{\tau}_{n,\text{OL}}^2(\lambda_m))^2$. We then select the value of $_s\lambda_m$, say $\widehat{_s\lambda_m^{\text{opt}}}$, which minimizes this data-based MSE and take

$$\widehat{_s\lambda_n^{\text{opt}}} = \widehat{_s\lambda_m^{\text{opt}}}\{|R_n|/|\lambda_m R_0|\}^{1/(2+d)}.$$

We use Theorem 5.1 to appropriately recalibrate an estimate $\widehat{_s\lambda_m^{\text{opt}}}$ to estimate optimal subsample scaling for $R_n$-size regions. For optimal scaling estimation with NOL subsamples, we replace $\hat{\tau}_{n,\text{OL}}^2(\lambda_m), \hat{\tau}_{i,m,\text{OL}}^2$ with $\hat{\tau}_{n,\text{NOL}}^2(\lambda_m), \hat{\tau}_{i,m,\text{NOL}}^2$ above. Garcia-Soidan and Hall (1997) apply a similar empirical MSE selection procedure with subsample-based distribution estimators on rectangular sampling regions in $\mathbb{R}^2$.

An advantage of a plug-in estimate of scaling is that it is computationally less demanding than minimization of an estimated MSE. A nonparametric plug-in (NPI) procedure involves substituting estimates of unknown r.f. parameters appearing in $_s\lambda_n^{\text{opt}}$ from Theorem 5.1. To do this, we propose using subsample variance estimators based on two smoothing parameter choices. Let $\hat{\tau}_n^2(_s\lambda_n)$ denote a subsample variance estimator with scaling $_s\lambda_n$ in (4.1). Using a pilot scalar $_s\lambda_n^{(1)} = c_1 |R_n|^{1/(d+2)}$, $c_1 > 0$, we estimate the limiting variance $\tau^2$ appearing in $_s\lambda_n^{\text{opt}}$ with $\hat{\tau}_n^2(_s\lambda_n^{(1)})$. With a second smoothing parameter $_s\lambda_n^{(2)} = c_2 |R_n|^{1/(d+4)}$, $c_2 > 0$, we estimate the bias component $B_0$ with $\hat{B}_0 = 2_s\lambda_n^{(2)}[\hat{\tau}_n^2(2_s\lambda_n^{(2)}) - \hat{\tau}_n^2(_s\lambda_n^{(2)})]$. It follows easily from Theorems 3.1–4.1 that the estimator $\hat{B}_0$ is consistent when the bias of $\hat{\tau}_n^2(_s\lambda_n)$ is $-_s\lambda_n^{-1} B_0(1 + o(1))$. With time series $d = 1$, Lahiri, Furukawa and Lee (2003) suggest a similar bias estimate for the MBB variance estimator and show the order of $_s\lambda_n^{(2)}$ above is asymptotically optimal. Politis and Romano (1995) also consider combining two subsample estimators in kernel spectral density estimation. We conjecture that the order $_s\lambda_n^{(2)}$ is optimal for minimizing the asymptotic MSE in estimating $B_0$ with spatial subsampling ($d \geq 2$) and this can be established for rectangular sampling regions.

For subsample variance estimation of a time series mean, other plug-in rules for $_s\lambda_n^{\mathrm{opt}}$ are given in Carlstein (1986) [with AR(1) models], Léger, Politis and Romano (1992) and Politis and Romano (1993b).

## 8. Numerical studies.

8.1. *Performance comparison of subsample types.* We conducted a simulation study to compare the finite sample performances of OL and NOL subsample variance estimators of $\tau_n^2 = N_n \operatorname{Var}(\hat{\theta}_n)$, where $\hat{\theta}_n = \bar{Z}_{N_n}$ is the real-valued sample mean over a sampling region $R_n \subset \mathbb{R}^2$. Rectangular and circular regions $R_n$ of two different sizes were considered:

$$R_n := (-7, 7] \times (-9, 9], \qquad R_n := (-15, 15] \times (-21, 21],$$
$$R_n := \{\mathbf{x} \in \mathbb{R}^2 : \|\mathbf{x}\| \leq 9\}, \qquad R_n := \{\mathbf{x} \in \mathbb{R}^2 : \|\mathbf{x}\| \leq 20\}.$$

The smaller (larger) circle contains one $\mathbb{Z}^2$ integer point more (seven less) than the smaller (larger) rectangle. The rectangular regions have approximately the same ratio of side lengths.

Using the algorithm of Chan and Wood (1997), we generated mean zero Gaussian random fields on $\mathbb{Z}^2$ with one of the following covariance structures:

$$\text{Model } \mathrm{E}(\beta_1, \beta_2) : \sigma(\mathbf{k}) = \exp[-\beta_1|k_1| - \beta_2|k_2|],$$
(8.1) $$\quad \text{Model } \mathrm{G}(\beta_1, \beta_2) : \sigma(\mathbf{k}) = \exp[-\beta_1|k_1|^2 - \beta_2|k_2|^2],$$
$$\mathbf{k} = (k_1, k_2)' \in \mathbb{Z}^2, \beta_1, \beta_2 > 0.$$

Models E and G correspond to exponential and Gaussian covariograms, respectively. We consider the values $(\beta_1, \beta_2) = (0.5, 0.3), (1, 1)$ in both models to obtain isotropic and anisotropic covariograms exhibiting various rates of decay.

For each $R_n$ and covariance structure, we considered various amounts of subsample scaling $_s\lambda_n$ in the estimator $\hat{\tau}_n^2 \equiv \hat{\tau}_n^2(_s\lambda_n)$ based on OL or NOL subsamples. Here rectangular and circular subsamples correspond to translates of $_s\lambda_n R_0$ for $R_0 = (-1/2, 1/2]^2$, $\{\mathbf{x} \in \mathbb{R}^2 : \|\mathbf{x}\| \leq 1/2\}$. We estimated the normalized MSE, $\mathrm{E}(\hat{\tau}_n^2/\tau_n^2 - 1)^2$, listing results in Table 6 for Model E. (To save space, we omit similar tables for Model G, where the performance of the estimators was better.) Estimates of optimal scaling appear in Table 7. From these simulation results, we make the following observations:

TABLE 6

Normalized MSE $\mathrm{E}(\hat{\tau}_n^2/\tau_n^2 - 1)^2$ for OL/NOL subsample variance estimators $\hat{\tau}_n^2$ of $\tau_n^2 = N_n \mathrm{Var}(\bar{Z}_{N_n})$ on $R_n \cap \mathbb{Z}^2$ (based on 10,000 simulations). An asterisk (*) denotes a minimal MSE

| $_s\lambda_n$ | E(0.5, 0.3) | | E(1, 1) | | E(0.5, 0.3) | | E(1, 1) | |
|---|---|---|---|---|---|---|---|---|
| | OL | NOL | OL | NOL | OL | NOL | OL | NOL |
| | $R_n = (-7, 7] \times (-9, 9]$ | | | | $R_n = \{\mathbf{x} \in \mathbb{R}^2 : \|\mathbf{x}\| \leq 9\}$ | | | |
| 1 | 0.9074 | 0.9074 | 0.5855 | 0.5855 | 0.9075 | 0.9075 | 0.5871 | 0.5871 |
| 2 | 0.7645 | 0.7619 | 0.3312 | 0.3298 | 0.7413 | 0.7417 | 0.3303 | 0.3330 |
| 3 | 0.6367 | 0.6343 | 0.2201 | 0.2264 | 0.6386 | 0.6378 | 0.2252 | 0.2346* |
| 4 | 0.5490 | 0.5470 | 0.1926* | 0.2191* | 0.5991 | 0.6177 | 0.2332 | 0.2897 |
| 5 | 0.5051 | 0.5344 | 0.2106 | 0.3071 | 0.5255 | 0.5627 | 0.2126* | 0.3444 |
| 6 | 0.4999* | 0.4605* | 0.2533 | 0.2911 | 0.5246* | 0.4978* | 0.2567 | 0.3369 |
| 7 | 0.5242 | 0.4957 | 0.3086 | 0.4004 | 0.5311 | | | 0.2925 |
| | $R_n = (-15, 15] \times (-21, 21]$ | | | | $R_n = \{\mathbf{x} \in \mathbb{R}^2 : \|\mathbf{x}\| \leq 20\}$ | | | |
| 4 | 0.5290 | 0.5285 | 0.1820 | 0.1851 | 0.5849 | 0.5846 | 0.1825 | 0.1866 |
| 5 | 0.4370 | 0.4329 | 0.1170 | 0.1232 | 0.4743 | 0.4785 | 0.1186 | 0.1332 |
| 6 | 0.3693 | 0.3601 | 0.1115 | 0.1380 | 0.4180 | 0.4236 | 0.1119 | 0.1358 |
| 7 | 0.3226 | 0.3132 | 0.0983* | 0.1172* | 0.3698 | 0.3716 | 0.1007* | 0.1257* |
| 8 | 0.2931 | 0.2963 | 0.1061 | 0.1453 | 0.3313 | 0.3466 | 0.1055 | 0.1596 |
| 9 | 0.2777 | 0.2822 | 0.1085 | 0.1613 | 0.2901 | 0.3333 | 0.1119 | 0.2080 |
| 10 | 0.2734* | 0.2542* | 0.1298 | 0.2247 | 0.2849 | 0.3084* | 0.1254 | 0.2049 |
| 11 | 0.2779 | 0.3454 | 0.1388 | 0.2824 | 0.2803* | 0.3814 | 0.1397 | 0.3335 |
| 12 | 0.2891 | 0.3298 | 0.1680 | 0.2889 | 0.2868 | 0.3662 | 0.1596 | 0.3359 |

1. At optimal scaling, the MSEs of OL and NOL subsamples were similar. Under the strongest r.f. dependence in Model E(0.5, 0.3), NOL subsamples performed better. For the other covariogram models entailing weaker dependence, OL subsamples were always better.
2. Unlike with OL subsamples, the MSEs with NOL subsamples increased more rapidly when optimal scaling was not used. This implies estimation of $_s\lambda_n^{\mathrm{opt}}$ with OL subsamples is preferable.

TABLE 7

Optimal subsample scaling $_s\lambda_n^{\mathrm{opt}}$ for variance estimation of sample mean $\sqrt{N_n}\bar{Z}_{N_n}$ (determined from 10,000 simulations)

| $R_n$ | E(0.5, 0.3) | | G(0.5, 0.3) | | E(1, 1) | | G(1, 1) | |
|---|---|---|---|---|---|---|---|---|
| | OL | NOL | OL | NOL | OL | NOL | OL | NOL |
| $(-7, 7] \times (-9, 9]$ | 6 | 6 | 4 | 4 | 4 | 4 | 3 | 3 |
| $(-15, 15] \times (-21, 21]$ | 10 | 10 | 7 | 6 | 7 | 6 | 5 | 5 |
| $\{\mathbf{x} \in \mathbb{R}^2 : \|\mathbf{x}\| \leq 9\}$ | 6 | 6 | 5 | 3 | 5 | 3 | 3 | 3 |
| $\{\mathbf{x} \in \mathbb{R}^2 : \|\mathbf{x}\| \leq 20\}$ | 11 | 10 | 7 | 7 | 7 | 7 | 5 | 5 |

3. Table 7 shows that OL and NOL optimal scaling tended to be the same. NOL subsample scaling becomes clearly smaller in larger sample sizes; see also Table 4.
4. Optimal subsample scaling also decreased as the r.f. dependence structure weakened (e.g., faster decay of covariogram). In this case, the performance of the variance estimators also improved.

8.2. *Comparison of scaling estimation methods.* We also compared NPI and HJ estimation methods for scaling $_s\lambda_{n,\mathrm{OL}}^{\mathrm{opt}}$ with OL subsamples, using the covariogram models and sampling regions $R_n$ from Section 8.1. We again took the sample mean $\hat{\theta}_n = \bar{Z}_{N_n}$. For the NPI method, we chose smoothing parameters $c_1, c_2 \in \{0.5, 1, 2\}$. For each $R_n$, we used two pilot subsample sizes $\lambda_m$ for the HJ method. As a measure of performance of the NPI and HJ procedures, we considered the following quantity:

$$(8.2) \qquad \phi_n = \frac{\hat{\tau}_{n,\mathrm{OL}}^2(\widehat{_s\lambda_{n,\mathrm{OL}}^{\mathrm{opt}}}) - \hat{\tau}_{n,\mathrm{OL}}^2(_s\lambda_{n,\mathrm{OL}}^{\mathrm{opt}})}{\tau_n^2},$$

where $\hat{\tau}_{n,\mathrm{OL}}^2(_s\lambda_n)$ denotes the OL subsample variance estimator using scaling $_s\lambda_n$, $\widehat{_s\lambda_{n,\mathrm{OL}}^{\mathrm{opt}}}$ represents an estimate of optimal scaling $_s\lambda_{n,\mathrm{OL}}^{\mathrm{opt}}$, and $\tau_n^2$ is the variance parameter. Hence, $\phi_n$ measures the relative deviation of an OL subsample estimator of $\tau_n^2$ based on estimated scaling compared to the "best" OL subsample estimator. Values of $\phi_n$ near zero would suggest that $\hat{\tau}_{n,\mathrm{OL}}^2(\widehat{_s\lambda_{n,\mathrm{OL}}^{\mathrm{opt}}})$ performed nearly as well as the optimal subsample estimator $\hat{\tau}_{n,\mathrm{OL}}^2(_s\lambda_n^{\mathrm{opt}})$.

From the results reported partially in Table 8, the choices of smoothing parameters

$$c_2 = 0.5 \quad \text{and} \quad c_1 \in \{0.5, 1\}$$

gave good results for estimating $_s\lambda_n^{\mathrm{opt}}$ in the NPI approach. We recommend these values for implementing the NPI method. The HJ method also tended to perform better with smaller smoothing parameter choices $\lambda_m$, which agrees with the $\lambda_m$ selections of Hall and Jing (1996) for time series. (We chose $\lambda_m$ so that an estimated MSE could be maximized over at least five different $_s\lambda_m$ arguments.) Table 9 gives frequency distributions of estimated optimal scaling $_s\lambda_{n,\mathrm{OL}}^{\mathrm{opt}}$ under other covariogram models and regions $R_n$. Table 7 lists values of $_s\lambda_{n,\mathrm{OL}}^{\mathrm{opt}}$. These results indicate that the NPI and HJ procedures exhibit good finite sample properties in estimating $_s\lambda_n^{\mathrm{opt}}$ and are competitive.

Table 8

Values of $\mathrm{E}(\phi_n^2)$ for NPI and HJ methods (each based on 1000 simulations), where $\phi_n$ is as in (8.2). HJ method uses $(\lambda_{m_1}, \lambda_{m_2}) = (5,10), (7,14), (3,6), (4,8)$, respectively, on regions $R_n$ from left to right. Minimal MSE is denoted with an asterisk "*" for each $R_n$ and covariogram model

| $R_n$ | | $(-7,7] \times (-9,9]$ | | $(-15,15] \times (-21,21]$ | | $\{\mathbf{x} \in \mathbb{R}^2 : \|\mathbf{x}\| \leq 9\}$ | | $\{\mathbf{x} \in \mathbb{R}^2 : \|\mathbf{x}\| \leq 20\}$ | |
| $c_1$ | $c_2$ | E(0.5, 0.3) | G(1,1) | E(0.5, 0.3) | G(1,1) | E(0.5, 0.3) | G(1,1) | E(0.5, 0.3) | G(1,1) |
|---|---|---|---|---|---|---|---|---|---|
| 0.5 | 0.5 | 0.0022* | 0.0106 | 0.0025 | 0.0075 | 0.0013 | 0.0093 | 0.0015 | 0.0034 |
|  | 1 | 0.0654 | 0.0614 | 0.0296 | 0.0288 | 0.0405 | 0.0559 | 0.0139 | 0.0266 |
|  | 2 | 0.0703 | 0.2470 | 0.1044 | 0.1000 | 0.0405 | 0.2532 | 0.1628 | 0.0862 |
| 1 | 0.5 | 0.0299 | 0.0031* | 0.0101 | 0.0027 | 0.0118 | 0.0047* | 0.0334 | 0.0011* |
|  | 1 | 0.0065 | 0.0706 | 0.0019* | 0.0206 | 0.0030 | 0.0644 | 0.0006* | 0.0192 |
|  | 2 | 0.0412 | 0.2040 | 0.0317 | 0.0968 | 0.0233 | 0.2098 | 0.0600 | 0.0911 |
| 2 | 0.5 | 0.0412 | 0.0352 | 0.0369 | 0.0055 | 0.0212 | 0.0205 | 0.0709 | 0.0029 |
|  | 1 | 0.0040 | 0.1081 | 0.0051 | 0.0157 | 0.0010 | 0.0961 | 0.0152 | 0.0133 |
|  | 2 | 0.0439 | 0.2582 | 0.0278 | 0.1346 | 0.0255 | 0.2676 | 0.0134 | 0.1206 |
| HJ, $\lambda_{m_1}$ | | 0.0100 | 0.0098 | 0.0161 | 0.0001* | 0.0001* | 0.0709 | 0.0334 | 0.0288 |
| HJ, $\lambda_{m_2}$ | | 0.0178 | 0.1766 | 0.0048 | 0.0130 | 0.0069 | 0.0360 | 0.0630 | 0.0337 |

Table 9

Frequency distribution of estimated optimal OL subsample scaling with NPI and HJ methods (based on 1000 simulations). Along with $c_2 = 0.5$, NPI1 and NPI2 use $c_1 = 0.5$ and 1, respectively. True optimal scaling values ${}_s\lambda_{n,\mathrm{OL}}^{\mathrm{opt}}$ are given in Table 7

| | | Estimates $\widehat{{}_s\lambda_{n,\mathrm{OL}}^{\mathrm{opt}}}$ of optimal scaling ${}_s\lambda_{n,\mathrm{OL}}^{\mathrm{opt}}$ | | | | | | | | |
| $R_n$/Model | Method | 2 | 3 | 4 | 5 | 6 | 7 | 8 | 9 | 10 |
|---|---|---|---|---|---|---|---|---|---|---|
| $(-7,7] \times (-9,9]$ | NPI1 | | 98 | 901 | 1 | | | | | |
| E(1,1) | NPI2 | | 307 | 686 | 7 | | | | | |
| | HJ, $\lambda_m = 5$ | 150 | | 850 | | | | | | |
| $\{\mathbf{x} \in \mathbb{R}^2 : \|\mathbf{x}\| \leq 9\}$ | NPI1 | | | 7 | 993 | | | | | |
| G(0.5, 0.3) | NPI2 | | | 876 | 124 | | | | | |
| | HJ, $\lambda_m = 3$ | | | | | 963 | | 37 | | |
| $(-15,15] \times (-21,21]$ | NPI1 | | | 2 | 9 | 62 | 276 | 450 | 192 | 9 |
| E(1,1) | NPI2 | | | 1 | 14 | 241 | 726 | 18 | | |
| | HJ, $\lambda_m = 7$ | 1 | | | 856 | | 143 | | | |
| $\{\mathbf{x} \in \mathbb{R}^2 : \|\mathbf{x}\| \leq 20\}$ | NPI1 | | | | 2 | 21 | 272 | 590 | 115 | |
| G(0.5, 0.3) | NPI2 | | | | 2 | 134 | 864 | | | |
| | HJ, $\lambda_m = 4$ | | | | | | 723 | 277 | | |

**9. Proofs for variance expansions.** For the proofs, we use $\mathcal{C}$ to denote generic positive constants that do not depend on $n$ or any $\mathbb{Z}^d$ integers (or $\mathbf{Z}^d$ lattice points). The real number $r$, appearing in some proofs, always assumes the value stated under Condition $M_r$ with respect to the lemma or theorem under consideration. Unless otherwise specified, limits in order symbols are taken letting $n$ tend to infinity.

In the following, we denote the indicator function as $I_{\{\cdot\}}$ (i.e., $I_{\{\cdot\}} \in \{0,1\}$ and $I_{\{A\}} = 1$ if and only if an event $A$ holds). For two sequences $\{s_n\}$ and $\{t_n\}$ of positive real numbers, we write $s_n \sim t_n$ if $s_n/t_n \to 1$ as $n \to \infty$. We write $\lambda_n^{\max}$ and $_s\lambda_n^{\max}$ for the largest diagonal entries of $\Delta_n$ and $_s\Delta_n$, respectively, while $_s\lambda_{\min}^{(n)} \geq 1$ will denote the smallest diagonal entry of $_s\Delta_n$.

We require a few lemmas for the proofs.

LEMMA 9.1. *Suppose $T_1, T_2 \subset \mathbf{Z}^d \equiv \mathbf{t} + \mathbb{Z}^d$ are bounded. Let $p, q > 0$ where $1/p + 1/q < 1$. If $X_1$, $X_2$ are random variables, with $X_i$ measurable with respect to $\mathcal{F}_Z(T_i), i = 1, 2$, then*

$$|\operatorname{Cov}(X_1, X_2)| \leq 8(\mathrm{E}|X_1|^p)^{1/p}(\mathrm{E}|X_2|^q)^{1/q} \alpha\bigg(\operatorname{dis}(T_1, T_2); \max_{i=1,2}|T_i|\bigg)^{1-1/p-1/q},$$

*provided expectations are finite and $\operatorname{dis}(T_1, T_2) > 0$.*

The proof follows from Theorem 3, Doukhan [(1994), page 9].

LEMMA 9.2. *Let $r \in \mathbb{Z}_+$. Under Assumption A.3 and Condition $M_r$, for $1 \leq m \leq 2r$ and any $T \subset \mathbf{Z}^d \equiv \mathbf{t} + \mathbb{Z}^d$,*

$$\mathrm{E}\bigg\|\sum_{\mathbf{s} \in T}(Z(\mathbf{s}) - \mu)\bigg\|^m \leq \mathcal{C}(\alpha)|T|^{m/2};$$

$\mathcal{C}(\alpha)$ *is a constant that depends only on the coefficients $\alpha(k,l)$, $l \leq 2r$, and $\mathrm{E}\|Z(\mathbf{t})\|^{2r+\delta}$.*

The proof follows from Theorem 1, Doukhan [(1994), pages 26–31] and Jensen's inequality.

We next determine the asymptotic sizes of important sets relevant to the sampling or subsampling designs.

LEMMA 9.3. *Under Assumptions A.1 and A.2, the number of sampling sites within:*

(a) *the sampling region $R_n$: $N_n = |R_n \cap \mathbf{Z}^d| \sim |R_0| \cdot \det(\Delta_n)$;*
(b) *an OL subsample, $R_{\mathbf{i},n}$, $\mathbf{i} \in J_{\mathrm{OL}}$: $_sN_n \sim |R_0| \cdot \det(_s\Delta_n)$;*
(c) *a NOL subsample, $\tilde{R}_{\mathbf{i},n}$, $\mathbf{i} \in J_{\mathrm{NOL}}$: $_sN_{\mathbf{i},n} \sim |R_0| \cdot \det(_s\Delta_n)$.*

*The number of:*

(d) *OL subsamples within $R_n$: $|J_{\mathrm{OL}}| \sim |R_0| \cdot \det(\Delta_n)$;*
(e) *NOL subsamples within $R_n$: $|J_{\mathrm{NOL}}| \sim |R_0| \cdot \det(\Delta_n) \cdot \det(_s\Delta_n)^{-1}$;*

(f) *sampling sites near the border of a subsample, $R_{\mathbf{i},n}$ or $\tilde{R}_{\mathbf{i},n}$, is less than*

$$\sup_{\mathbf{i}\in\mathbb{Z}^d}|\{\mathbf{j}\in\mathbb{Z}^d:T^{\mathbf{j}}\cap\overline{R_{\mathbf{i},n}}\neq\varnothing,\ T^{\mathbf{j}}\cap\overline{R_{\mathbf{i},n}^c}\neq\varnothing\ for\ T^{\mathbf{j}}=\mathbf{j}+[-2,2]^d\}|$$
$$\leq \mathcal{C}(_s\lambda_n^{\max})^{d-1}.$$

Results follow from the boundary condition on $R_0$; see Nordman (2002) for more details.

We require the next lemma for counting the number of subsampling regions which are separated by an appropriately "small" integer translate; we shall apply this lemma in the proof of Theorem 3.1. For $\mathbf{k} = (k_1, \ldots, k_d)' \in \mathbb{Z}^d$, define the following sets:

$$J_n(\mathbf{k}) = |\{\mathbf{i} \in J_{\mathrm{OL}} : \mathbf{i} + \mathbf{k} + {}_s\Delta_n R_0 \subset R_n\}|,$$
$$E_n = \{\mathbf{k} \in \mathbb{Z}^d : |k_j| \leq {}_s\lambda_j^{(n)},\ j = 1, \ldots, d\}.$$

LEMMA 9.4. *Under Assumption A.2,*
$$\max_{\mathbf{k}\in E_n}\left|1 - \frac{J_n(\mathbf{k})}{|J_{\mathrm{OL}}|}\right| = o(1).$$

PROOF. For $\mathbf{k} \in E_n$, write the set $J_n^*(\mathbf{k})$ and bound its cardinality
$$J_n^*(\mathbf{k}) = |\{\mathbf{i} \in J_{\mathrm{OL}} : (\mathbf{i} + \mathbf{k} + {}_s\Delta_n R_0) \cap \Delta_n R_0^c \neq \varnothing\}|$$
$$\leq |\{\mathbf{i}\in\mathbb{Z}^d:T^{\mathbf{i}}\cap\Delta_n\overline{R_0^c}\neq\varnothing,\ T^{\mathbf{i}}\cap\Delta_n\overline{R_0}\neq\varnothing;\ T^{\mathbf{i}}=\mathbf{i}+{}_s\lambda_n^{\max}[-2,2]^d\}|$$
$$\leq \mathcal{C}_s\lambda_n^{\max}(\lambda_n^{\max})^{d-1}$$

by the boundary condition on $R_0$. We have then that for all $\mathbf{k} \in E_n$,
$$|J_{\mathrm{OL}}| \geq J_n(\mathbf{k}) = |J_{\mathrm{OL}}| - J_n^*(\mathbf{k}) \geq |J_{\mathrm{OL}}| - \mathcal{C}_s\lambda_n^{\max}(\lambda_n^{\max})^{d-1}.$$

By Assumption A.2 and the growth rate of $|J_{\mathrm{OL}}|$ from Lemma 9.3, the proof is complete. $\square$

We now provide a theorem which captures the main contribution to the asymptotic variance expansion of the OL subsample variance estimator $\hat{\tau}_{n,\mathrm{OL}}^2$ from Theorem 3.1.

THEOREM 9.1. *For $\mathbf{i} \in \mathbb{Z}^d$, let $Y_{\mathbf{i},n} = \nabla'(Z_{\mathbf{i},n} - \mu)$. Under the assumptions and conditions of Theorem* 3.1
$$_sN_n \sum_{\mathbf{k}\in E_n} \mathrm{Cov}(Y_{\mathbf{0},n}^2, Y_{\mathbf{k},n}^2) = K_0 \cdot [2\tau^4](1 + o(1)),$$
*where the constant $K_0$ is defined in Theorem* 3.1.

PROOF. We give only a sketch of the important features; for more details, see Nordman (2002). For a set $T \subset \mathbb{R}^d$, define the function $\overline{\Sigma}(\cdot)$ as

$$\overline{\Sigma}(T) = \sum_{\mathbf{s} \in \mathbf{Z}^d \cap T} \nabla'(Z(\mathbf{s}) - \mu).$$

With the set intersection $R_{\mathbf{k},n}^{(\mathrm{I})} = {}_sR_n \cap (\mathbf{k} + {}_sR_n)$, $\mathbf{k} \in \mathbb{Z}^d$, write functions

$$H_{1n}(\mathbf{k}) = \overline{\Sigma}(R_{\mathbf{k},n} \setminus R_{\mathbf{k},n}^{(\mathrm{I})}),$$

$$H_{2n}(\mathbf{k}) = \overline{\Sigma}(R_{\mathbf{0},n} \setminus R_{\mathbf{k},n}^{(\mathrm{I})}),$$

$$H_{3n}(\mathbf{k}) = \overline{\Sigma}(R_{\mathbf{k},n}^{(\mathrm{I})}).$$

These represent, respectively, sums over sites in $R_{\mathbf{k},n}$ but not $R_{\mathbf{0},n} = {}_sR_n$, $R_{\mathbf{0},n}$ but not $R_{\mathbf{k},n}$ and both $R_{\mathbf{0},n}$ and $R_{\mathbf{k},n}$. Then define $h_n(\cdot) : \mathbb{Z}^d \to \mathbb{R}$ as

$$h_n(\mathbf{k}) = \mathrm{E}[H_{1n}^2(\mathbf{k})]\mathrm{E}[H_{2n}^2(\mathbf{k})] + \mathrm{E}[H_{1n}^2(\mathbf{k})]\mathrm{E}[H_{3n}^2(\mathbf{k})]$$
$$+ \mathrm{E}[H_{2n}^2(\mathbf{k})]\mathrm{E}[H_{3n}^2(\mathbf{k})] + \mathrm{E}[H_{3n}^4(\mathbf{k})] - ({}_sN_n)^4[\mathrm{E}(Y_{\mathbf{0},n}^2)]^2.$$

We will make use of the following proposition.

PROPOSITION 9.1. *Under the assumptions and conditions of Theorem* 3.1,

$$\max_{\mathbf{k} \in E_n} |({}_sN_n)^2 \mathrm{Cov}(Y_{\mathbf{0},n}^2, Y_{\mathbf{k},n}^2) - ({}_sN_n)^{-2} h_n(\mathbf{k})| = o(1).$$

The proof of Proposition 9.1 can be found in Nordman (2002) and involves cutting out $\mathbf{Z}^d$ lattice points near the borders of $R_{\mathbf{0},n}$ and $R_{\mathbf{k},n}$, say, $B_{\mathbf{0},n}$ and $B_{\mathbf{k},n}$ with

(9.1)
$$\ell_n = \lfloor ({}_s\lambda_{\min}^{(n)})^e \rfloor,$$
$$B_{\mathbf{j},n} = \{\mathbf{i} \in \mathbf{Z}^d : \mathbf{i} \in \overline{R_{\mathbf{j},n}}, (\mathbf{i} + \ell_n(-1,1]^d) \cap \overline{R_{\mathbf{j},n}^c} \neq \varnothing\}, \qquad \mathbf{j} \in \mathbb{Z}^d,$$

where $e = (\kappa\delta/\{(2r+\delta)(2r-1-1/d)\}+1)/2 < 1$ from Condition $M_r$. Here $\ell_n \to \infty$, $\ell_n = o({}_s\lambda_{\min}^{(n)})$ is chosen so that the remaining observations in $R_{\mathbf{0},n}, R_{\mathbf{k},n}, R_{\mathbf{k},n}^{(\mathrm{I})}$ are nearly independent upon removing $B_{\mathbf{0},n}, B_{\mathbf{k},n}$ points and, using the $R_0$-boundary condition, the set cardinalities $|B_{\mathbf{0},n}|, |B_{\mathbf{k},n}| \leq C\ell_n({}_s\lambda_n^{\max})^{d-1}$ are of smaller order than ${}_sN_n$ (namely, these sets are asymptotically negligible in size).

By Proposition 9.1 and $|E_n| = O({}_sN_n)$, we have

(9.2) $$\left| {}_sN_n \sum_{\mathbf{k} \in E_n} \mathrm{Cov}(Y_{\mathbf{0},n}^2, Y_{\mathbf{k},n}^2) - ({}_sN_n)^{-3} \sum_{\mathbf{k} \in E_n} h_n(\mathbf{k}) \right| = o(1).$$

Consequently, we need only focus on $(_sN_n)^{-3}\sum_{\mathbf{k}\in E_n}h_n(\mathbf{k})$ to complete the proof of Theorem 9.1.

For measurability reasons, we create a set defined in terms of the $\mathbb{R}^d$ Lebesgue measure,

$$\mathrm{E}^+ \equiv (0,1) \cap \left\{\varepsilon < \frac{\det(\Delta_0)|R_0|}{2}:\right.$$

$$\left.|\{\mathbf{x}\in\mathbb{R}^d:|(\mathbf{x}+\Delta_0 R_0)\cap \Delta_0 R_0|=\varepsilon \text{ or } \det(\Delta_0)|R_0|-\varepsilon\}|=0\right\}.$$

Note the set $(0,1)\cap(0,\det(\Delta_0)|R_0|/2)\setminus \mathrm{E}^+$ is at most countable [cf. Billingsley (1986), Theorem 10.4]. For $\varepsilon \in \mathrm{E}^+$, define a new set as a function of $\varepsilon$ and $n$:

$$\tilde{R}_{\varepsilon,n} = \{\mathbf{k}\in\mathbb{Z}^d: |R^{(\mathrm{I})}_{\mathbf{k},n}| > \varepsilon(_s\lambda_1^{(n)})^d,\ |_sR_n\setminus R^{(\mathrm{I})}_{\mathbf{k},n}| > \varepsilon(_s\lambda_1^{(n)})^d\}.$$

Here $\tilde{R}_{\varepsilon,n} \subset E_n$ because $\mathbf{k}\notin E_n$ implies $R^{(\mathrm{I})}_{\mathbf{k},n} = \varnothing$.

We now further simplify $(_sN_n)^{-3}\sum_{\mathbf{k}\in E_n}h_n(\mathbf{k})$ using the following proposition involving $\tilde{R}_{\varepsilon,n}$.

PROPOSITION 9.2. *There exist* $\mathrm{N}\in\mathbb{Z}_+$ *and a function* $b(\cdot):\mathrm{E}^+ \to (0,\infty)$ *such that* $b(\varepsilon)\downarrow 0$ *as* $\varepsilon\downarrow 0$ *and*

$$(9.3)\quad (_sN_n)^{-3}\left|\sum_{\mathbf{k}\in E_n}h_n(\mathbf{k}) - \sum_{\mathbf{k}\in \tilde{R}_{\varepsilon,n}}h_n(\mathbf{k})\right| \leq \mathcal{C}(\varepsilon + (_s\lambda_1^{(n)})^{-1} + [b(\varepsilon)]^d),$$

*where* $\mathcal{C}>0$ *does not depend on* $\varepsilon\in\mathrm{E}^+$ *or* $n\geq \mathrm{N}$.

The proof of Proposition 9.2 is tedious and given in Nordman (2002). The argument involves bounding the sum of $h_n(\cdot)$ over two separate sets in $E_n$: those integers in $E_n$ that are either "too large" or "too small" in magnitude to be included in $\tilde{R}_{\varepsilon,n}$.

To finish the proof, our approach (for an arbitrary $\varepsilon\in\mathrm{E}^+$) will be to write $(_sN_n)^{-3}\sum_{\mathbf{k}\in\tilde{R}_{\varepsilon,n}}h_n(\mathbf{k})$ as an integral of a step function $f_{\varepsilon,n}(\mathbf{x})$ with respect to the Lebesgue measure, then show $\lim_{n\to\infty}f_{\varepsilon,n}(\mathbf{x})$ exists almost everywhere (a.e.) on $\mathbb{R}^d$, and apply the Lebesgue dominated convergence theorem (LDCT). By letting $\varepsilon\downarrow 0$, we will obtain the limit of $_sN_n\sum_{\mathbf{k}\in E_n}\mathrm{Cov}(Y^2_{\mathbf{0},n},Y^2_{\mathbf{k},n})$.

Fix $\varepsilon\in\mathrm{E}^+$. With counting arguments based on the boundary condition of $R_0$ and the definition of $\tilde{R}_{\varepsilon,n}$, it holds that for some $\mathrm{N}_\varepsilon\in\mathbb{Z}_+$ and all $\mathbf{k}\in\tilde{R}_{\varepsilon,n}$: $|R^{(\mathrm{I})}_{\mathbf{k},n}\cap\mathbf{Z}^d|\geq 1$ and $_sN_n - |R^{(\mathrm{I})}_{\mathbf{k},n}\cap\mathbf{Z}^d|\geq 1$ when $n\geq\mathrm{N}_\varepsilon$. We can rewrite $(_sN_n)^{-2}h_n(\mathbf{k})$, $\mathbf{k}\in\tilde{R}_{\varepsilon,n}$, in the well-defined form (for $n\geq\mathrm{N}_\varepsilon$)

$$\frac{h_n(\mathbf{k})}{(_sN_n)^2} = \mathrm{E}\left[\frac{H^2_{1n}(\mathbf{k})}{_sN_n - |R^{(\mathrm{I})}_{\mathbf{k},n}\cap\mathbf{Z}^d|}\right]\mathrm{E}\left[\frac{H^2_{2n}(\mathbf{k})}{_sN_n - |R^{(\mathrm{I})}_{\mathbf{k},n}\cap\mathbf{Z}^d|}\right]\left(1 - \frac{|R^{(\mathrm{I})}_{\mathbf{k},n}\cap\mathbf{Z}^d|}{_sN_n}\right)^2$$

$$+ \sum_{j=1}^{2} \mathrm{E}\left[\frac{H_{jn}^2(\mathbf{k})}{{}_sN_n - |R_{\mathbf{k},n}^{(\mathrm{I})} \cap \mathbf{Z}^d|}\right] \mathrm{E}\left[\frac{H_{3n}^2(\mathbf{k})}{|R_{\mathbf{k},n}^{(\mathrm{I})} \cap \mathbf{Z}^d|}\right]$$

$$\times \left(1 - \frac{|R_{\mathbf{k},n}^{(\mathrm{I})} \cap \mathbf{Z}^d|}{{}_sN_n}\right)\left(\frac{|R_{\mathbf{k},n}^{(\mathrm{I})} \cap \mathbf{Z}^d|}{{}_sN_n}\right)$$

$$+ \mathrm{E}\left[\frac{H_{3n}^4(\mathbf{k})}{|R_{\mathbf{k},n}^{(\mathrm{I})} \cap \mathbf{Z}^d|^2}\right]\left(\frac{|R_{\mathbf{k},n}^{(\mathrm{I})} \cap \mathbf{Z}^d|}{{}_sN_n}\right)^2 - [{}_sN_n \mathrm{E}(Y_{\mathbf{0},n}^2)]^2.$$

For $\mathbf{x} = (x_1, \ldots, x_d)' \in \mathbb{R}^d$, write $\lfloor \mathbf{x} \rfloor = (\lfloor x_1 \rfloor, \ldots, \lfloor x_d \rfloor)' \in \mathbb{Z}^d$ and $\mathbf{x}_n = \lfloor {}_s\lambda_1^{(n)} \mathbf{x} \rfloor$. Let $f_{\varepsilon,n}(\mathbf{x}) : \mathbb{R}^d \to \mathbb{R}$ be the step function defined as

$$f_{\varepsilon,n}(\mathbf{x}) = ({}_sN_n)^{-2} I_{\{\mathbf{x}_n \in \tilde{R}_{\varepsilon,n}\}} h_n(\mathbf{x}_n).$$

We have then that (with the same fixed $\varepsilon \in \mathrm{E}^+$)

$$(9.4) \qquad \frac{1}{{}_sN_n} \sum_{\mathbf{k} \in \tilde{R}_{\varepsilon,n}} ({}_sN_n)^{-2} h_n(\mathbf{k}) = \frac{({}_s\lambda_1^{(n)})^d}{{}_sN_n} \int_{\mathbb{R}^d} f_{\varepsilon,n}(\mathbf{x}) \, d\mathbf{x}.$$

We focus on showing

$$\lim_{n \to \infty} f_{\varepsilon,n}(\mathbf{x}) = f_{\varepsilon}(\mathbf{x})$$

(9.5)

$$\equiv I_{\{\mathbf{x} \in \tilde{R}_{\varepsilon}\}} [2\tau^4]\left(\frac{|(\mathbf{x} + \Delta_0 R_0) \cap \Delta_0 R_0|}{\det(\Delta_0)|R_0|}\right)^2 \qquad \text{a.e. } \mathbf{x} \in \mathbb{R}^d,$$

with $\tilde{R}_{\varepsilon} = \{\mathbf{x} \in \mathbb{R}^d : |(\mathbf{x} + \Delta_0 R_0) \cap \Delta_0 R_0| > \varepsilon, |\Delta_0 R_0 \setminus (\mathbf{x} + \Delta_0 R_0)| > \varepsilon\}$ a Borel measurable set.

To establish (9.5), we begin by showing convergence of indicator functions

$$(9.6) \qquad I_{\{\mathbf{x}_n \in \tilde{R}_{\varepsilon,n}\}} \to I_{\{\mathbf{x} \in \tilde{R}_{\varepsilon}\}} \qquad \text{a.e. } \mathbf{x} \in \mathbb{R}^d.$$

Define the sets $A_n(\mathbf{x}) = ({}_s\lambda_1^{(n)})^{-1}\{(\mathbf{x}_n + {}_sR_n) \cap {}_sR_n\}$, $\tilde{A}_n(\mathbf{x}) = \{({}_s\lambda_1^{(n)})^{-1} {}_sR_n\} \setminus A_n(\mathbf{x})$ as a function of $\mathbf{x} \in \mathbb{R}^d$. The LDCT can be applied to show that for each $\mathbf{x} \in \mathbb{R}^d$, $|A_n(\mathbf{x})| \to |(\mathbf{x} + \Delta_0 R_0) \cap \Delta_0 R_0|$ and $|\tilde{A}_n(\mathbf{x})| \to |\Delta_0 R_0 \setminus (\mathbf{x} + \Delta_0 R_0)|$. Thus, if $\mathbf{x} \in \tilde{R}_{\varepsilon}$, then

$$(9.7) \qquad \begin{aligned} |A_n(\mathbf{x})| &\to |(\mathbf{x} + \Delta_0 R_0) \cap \Delta_0 R_0| > \varepsilon, \\ |\tilde{A}_n(\mathbf{x})| &\to |\Delta_0 R_0 \setminus (\mathbf{x} + \Delta_0 R_0)| > \varepsilon, \end{aligned}$$

implying further that $1 = I_{\{\mathbf{x}_n \in \tilde{R}_{\varepsilon,n}\}} \to I_{\{\mathbf{x} \in \tilde{R}_{\varepsilon}\}} = 1$ as $n \to \infty$. Now consider $\tilde{R}_{\varepsilon}^c$. If $\mathbf{x} \notin \tilde{R}_{\varepsilon}$ such that $|(\mathbf{x} + \Delta_0 R_0) \cap \Delta_0 R_0| < \varepsilon$ [or $|\Delta_0 R_0 \setminus (\mathbf{x} + \Delta_0 R_0)| < \varepsilon$], then $|A_n(\mathbf{x})| < \varepsilon$ [or $|\tilde{A}_n(\mathbf{x})| < \varepsilon$] eventually for large $n$ and $0 = I_{\{\mathbf{x}_n \in \tilde{R}_{\varepsilon,n}\}} \to I_{\{\mathbf{x} \in \tilde{R}_{\varepsilon}\}} = 0$ in this case. Finally, $\varepsilon \in \mathrm{E}^+$ implies that a last possible subset of $\tilde{R}_\varepsilon^c$ has Lebesgue measure zero; namely, $|\{\mathbf{x} \in \tilde{R}_\varepsilon^c : |(\mathbf{x} + \Delta_0 R_0) \cap \Delta_0 R_0| = \varepsilon$ or $|\Delta_0 R_0 \setminus (\mathbf{x} + \Delta_0 R_0)| = \varepsilon\}| = 0$. We have now proven (9.6).

We next establish a limit for $(_sN_n)^{-2}h_n(\mathbf{x}_n)$, $\mathbf{x} \in \tilde{R}_\varepsilon$. We wish to show

$$(9.8) \qquad \frac{|R^{(\mathrm{I})}_{\mathbf{x}_n,n} \cap \mathbf{Z}^d|}{_sN_n} \to \frac{|(\mathbf{x} + \Delta_0 R_0) \cap \Delta_0 R_0|}{\det(\Delta_0)|R_0|}, \qquad \mathbf{x} \in \tilde{R}_\varepsilon.$$

Using the bound $||R^{(\mathrm{I})}_{\mathbf{x}_n,n}| - |R^{(\mathrm{I})}_{\mathbf{x}_n,n} \cap \mathbf{Z}^d|| \le \mathcal{C}(_s\lambda_n^{\max})^{d-1}$ from the $R_0$-boundary condition and noting the limit in (9.7) for $(_s\lambda_1^{(n)})^{-d}|R^{(\mathrm{I})}_{\mathbf{x}_n,n}| = |A_n(\mathbf{x})|$, we find $(_s\lambda_1^{(n)})^{-d}|R^{(\mathrm{I})}_{\mathbf{x}_n,n} \cap \mathbf{Z}^d| \to |(\mathbf{x} + \Delta_0 R_0) \cap \Delta_0 R_0|$, $\mathbf{x} \in \tilde{R}_\varepsilon$. By this and $(_s\lambda_1^{(n)})^d/_sN_n \to (\det(\Delta_0)|R_0|)^{-1}$, (9.8) follows.

We can also establish: for each $\mathbf{x} \in \tilde{R}_\varepsilon$, $j = 1$ or $2$,

$$(9.9) \qquad \begin{aligned} \mathrm{E}\left[\frac{H_{3n}^{2j}(\mathbf{x}_n)}{|R^{(\mathrm{I})}_{\mathbf{x}_n,n} \cap \mathbf{Z}^d|^j}\right] &\to \mathrm{E}([\nabla' Z_\infty]^{2j}), \\ \mathrm{E}\left[\frac{H_{jn}^2(\mathbf{x}_n)}{_sN_n - |R^{(\mathrm{I})}_{\mathbf{x}_n,n} \cap \mathbf{Z}^d|}\right], \qquad &_sN_n\mathrm{E}(Y_{\mathbf{0},n}^2) \to \mathrm{E}([\nabla' Z_\infty]^2), \end{aligned}$$

where $\nabla' Z_\infty$ is a normal $\mathcal{N}(0, \tau^2)$ random variable and so it follows that $\mathrm{E}([\nabla' Z_\infty]^{2j}) = (2j-1)\tau^{2j}$, $j = 1, 2$. The limits in (9.9) follow essentially from the central limit theorem (CLT) of Bolthausen (1982), after verifying that the CLT can be applied; see Nordman (2002) for more details.

Putting (9.6), (9.8) and (9.9) together, we have shown the (a.e.) convergence of the univariate functions $f_{\varepsilon,n}(\mathbf{x})$ as in (9.5). For $\mathbf{k} \in E_n$ and $n \ge \mathrm{N}_\varepsilon$, Lemma 9.2 ensures: $(_sN_n)^{-2}|h_n(\mathbf{k})| \le \mathcal{C}$, implying that for $\mathbf{x} \in \mathbb{R}^d : |f_{\varepsilon,n}(\mathbf{x})| \le \mathcal{C}I_{\{\mathbf{x} \in [-c,c]^d\}}$ for some $c > 0$ by Assumption A.1. With this uniform bound on $f_{\varepsilon,n}(\cdot)$ and the limits in (9.5), we can apply the LDCT to get

$$(9.10) \qquad \lim_{n\to\infty} \int_{\mathbb{R}^d} f_{\varepsilon,n}(\mathbf{x})\, d\mathbf{x} = \int_{\mathbb{R}^d} f_\varepsilon(\mathbf{x})\, d\mathbf{x}, \qquad \varepsilon \in \mathrm{E}^+.$$

Let $\{\varepsilon_m\}_{m=1}^\infty \subset \mathrm{E}^+$ where $\varepsilon_m \downarrow 0$. Then $\tilde{R}_{\varepsilon_m} \subset \Delta_0[-1,1]^d$ and $\lim_{m\to\infty} I_{\{\mathbf{x} \in \tilde{R}_{\varepsilon_m}\}} \to I_{\{\mathbf{x} \in \tilde{R}_0\}}$ for $\mathbf{x} \ne \mathbf{0} \in \mathbb{R}^d$, with $\tilde{R}_0 = \{\mathbf{x} \in \mathbb{R}^d : 0 < |(\mathbf{x} + \Delta_0 R_0) \cap \Delta_0 R_0| < \det(\Delta_0)|R_0|\}$. Hence, by the LDCT,

$$(9.11) \qquad \begin{aligned} \lim_{m\to\infty} \int_{\mathbb{R}^d} f_{\varepsilon_m}(\mathbf{x})\, d\mathbf{x} &= \int_{\mathbb{R}^d} f_0(\mathbf{x})\, d\mathbf{x}, \\ f_0(\mathbf{x}) &\equiv I_{\{\mathbf{x} \in \tilde{R}_0\}}[2\tau^4]\left(\frac{|(\mathbf{x} + \Delta_0 R_0) \cap \Delta_0 R_0|}{\det(\Delta_0)|R_0|}\right)^2. \end{aligned}$$

From (9.2)–(9.4), (9.10) and (9.11) and $({}_s\lambda_1^{(n)})^d/{}_sN_n \to (\det(\Delta_0)|R_0|)^{-1}$, we have that

$$\limsup_{n\to\infty}\left|{}_sN_n \sum_{\mathbf{k}\in E_n} \mathrm{Cov}(Y_{\mathbf{0},n}^2, Y_{\mathbf{k},n}^2) - \frac{1}{\det(\Delta_0)|R_0|}\int_{\mathbb{R}^d} f_0(\mathbf{x})\,d\mathbf{x}\right|$$

$$\leq \limsup_{n\to\infty}\left|{}_sN_n \sum_{\mathbf{k}\in E_n} \mathrm{Cov}(Y_{\mathbf{0},n}^2, Y_{\mathbf{k},n}^2) - \frac{({}_s\lambda_1^{(n)})^d}{{}_sN_n}\int_{\mathbb{R}^d} f_{\varepsilon_m,n}(\mathbf{x})\,d\mathbf{x}\right|$$

$$+ \frac{1}{\det(\Delta_0)|R_0|}\left|\int_{\mathbb{R}^d} f_{\varepsilon_m}(\mathbf{x}) - f_0(\mathbf{x})\,d\mathbf{x}\right|$$

$$+ \limsup_{n\to\infty}\left|\frac{({}_s\lambda_1^{(n)})^d}{{}_sN_n}\int_{\mathbb{R}^d} f_{\varepsilon_m,n}(\mathbf{x})\,d\mathbf{x} - \frac{1}{\det(\Delta_0)|R_0|}\int_{\mathbb{R}^d} f_{\varepsilon_m}(\mathbf{x})\,d\mathbf{x}\right|$$

$$\leq \mathcal{C}(\varepsilon_m + [b(\varepsilon_m)]^d) + \frac{1}{\det(\Delta_0)|R_0|}\left|\int_{\mathbb{R}^d} f_{\varepsilon_m}(\mathbf{x}) - f_0(\mathbf{x})\,d\mathbf{x}\right|$$

$$\to 0 \quad \text{as } \varepsilon_m \downarrow 0.$$

Finally,

$$\frac{1}{\det(\Delta_0)|R_0|}\int_{\mathbb{R}^d} f_0(\mathbf{x})\,d\mathbf{x} = \frac{2\tau^4}{|R_0|}\int_{\mathbb{R}^d} \frac{|(\mathbf{y}+R_0)\cap R_0|^2}{|R_0|^2}\,d\mathbf{y},$$

using a change of variables $\mathbf{y} = \Delta_0^{-1}\mathbf{x}$. This completes the proof of Theorem 9.1. □

For clarity of exposition, we will prove Theorem 3.1, parts (a) and (b), separately for the OL and NOL subsample variance estimators.

9.1. *Proof of Theorem 3.1(a).* For $\mathbf{i}\in J_{\mathrm{OL}}$, we use a Taylor expansion of $H(\cdot)$ (around $\mu$) to rewrite the statistic $\hat{\theta}_{\mathbf{i},n}^{\mathrm{OL}} = H(Z_{\mathbf{i},n})$,

$$\hat{\theta}_{\mathbf{i},n}^{\mathrm{OL}} = H(\mu) + \sum_{\|\alpha\|_1=1} c_\alpha (Z_{\mathbf{i},n}-\mu)^\alpha$$

(9.12)
$$+ 2\sum_{\|\alpha\|_1=2} \frac{(Z_{\mathbf{i},n}-\mu)^\alpha}{\alpha!}\int_0^1 (1-\omega)D^\alpha H(\mu+\omega(Z_{\mathbf{i},n}-\mu))\,d\omega$$

$$\equiv H(\mu) + Y_{\mathbf{i},n} + Q_{\mathbf{i},n}.$$

We also have

$$\tilde{\theta}_n^{\mathrm{OL}} = H(\mu) + |J_{\mathrm{OL}}|^{-1}\sum_{\mathbf{i}\in J_{\mathrm{OL}}} Y_{\mathbf{i},n} + |J_{\mathrm{OL}}|^{-1}\sum_{\mathbf{i}\in J_{\mathrm{OL}}} Q_{\mathbf{i},n} \equiv H(\mu) + \bar{Y}_n + \bar{Q}_n.$$

Then
$$\hat{\tau}_{n,\mathrm{OL}}^2 = {}_sN_n\left[\frac{1}{|J_{\mathrm{OL}}|}\sum_{\mathbf{i}\in J_{\mathrm{OL}}} Y_{\mathbf{i},n}^2 + \frac{1}{|J_{\mathrm{OL}}|}\sum_{\mathbf{i}\in J_{\mathrm{OL}}} Q_{\mathbf{i},n}^2 \right.$$
$$\left. + \frac{2}{|J_{\mathrm{OL}}|}\sum_{\mathbf{i}\in J_{\mathrm{OL}}} Y_{\mathbf{i},n}Q_{\mathbf{i},n} - \bar{Y}_n^2 - \bar{Q}_n^2 - 2(\bar{Y}_n)(\bar{Q}_n)\right].$$

We establish Theorem 3.1(a) in two parts by showing

(9.13)
$$\text{(a)} \quad \mathrm{Var}\left(\frac{{}_sN_n}{|J_{\mathrm{OL}}|}\sum_{\mathbf{i}\in J_{\mathrm{OL}}} Y_{\mathbf{i},n}^2\right) = K_0 \cdot \frac{\det({}_s\Delta_n)}{\det(\Delta_n)} \cdot [2\tau^4](1+o(1)),$$
$$\text{(b)} \quad \left|\mathrm{Var}(\hat{\tau}_{n,\mathrm{OL}}^2) - \mathrm{Var}\left(\frac{{}_sN_n}{|J_{\mathrm{OL}}|}\sum_{\mathbf{i}\in J_{\mathrm{OL}}} Y_{\mathbf{i},n}^2\right)\right| = o\left(\frac{\det({}_s\Delta_n)}{\det(\Delta_n)}\right).$$

We will begin with proving (9.13)(a). For $\mathbf{k}\in\mathbb{Z}^d$, let $\sigma_n(\mathbf{k}) = \mathrm{Cov}(Y_{\mathbf{0},n}^2, Y_{\mathbf{k},n}^2)$. We write
$$\frac{({}_sN_n)^2}{|J_{\mathrm{OL}}|^2}\mathrm{Var}\left(\sum_{\mathbf{i}\in J_{\mathrm{OL}}} Y_{\mathbf{i},n}^2\right) = \frac{({}_sN_n)^2}{|J_{\mathrm{OL}}|^2}\left(\sum_{\mathbf{k}\in E_n} J_n(\mathbf{k})\sigma_n(\mathbf{k}) + \sum_{\mathbf{k}\in\mathbb{Z}^d\setminus E_n} J_n(\mathbf{k})\sigma_n(\mathbf{k})\right)$$
$$\equiv W_{1n} + W_{2n}.$$

By stationarity and Lemma 9.2, we bound $|\sigma_n(\mathbf{k})| \leq \mathrm{E}(Y_{\mathbf{0},n}^4) \leq \mathcal{C}({}_sN_n)^{-2}$, $\mathbf{k}\in\mathbb{Z}^d$. Using this covariance bound, Lemmas 9.3 and 9.4 and $|E_n|\leq 3^d\det({}_s\Delta_n)$,

(9.14)
$$\left|\frac{({}_sN_n)^2}{|J_{\mathrm{OL}}|}\sum_{\mathbf{k}\in E_n}\sigma_n(\mathbf{k}) - W_{1n}\right| \leq \mathcal{C}\frac{|E_n|}{|J_{\mathrm{OL}}|}\cdot\max_{\mathbf{k}\in E_n}\left|1 - \frac{J_n(\mathbf{k})}{|J_{\mathrm{OL}}|}\right|$$
$$= o\left(\frac{\det({}_s\Delta_n)}{\det(\Delta_n)}\right).$$

Then applying Theorem 9.1 and Lemma 9.3,

(9.15)
$$\frac{({}_sN_n)^2}{|J_{\mathrm{OL}}|}\sum_{\mathbf{k}\in E_n}\sigma_n(\mathbf{k}) = K_0 \cdot \frac{\det({}_s\Delta_n)}{\det(\Delta_n)}\cdot[2\tau^4](1+o(1)).$$

By (9.14) and (9.15), we need only show that $W_{2n} = o\left(\det({}_s\Delta_n)/\det(\Delta_n)\right)$ to finish the proof of (9.13)(a).

For $\mathbf{i}\in\mathbb{Z}^d$, denote a set of lattice points within a translated rectangular region:
$$F_{\mathbf{i},n} = \left(\mathbf{i} + \prod_{j=1}^d(-\lceil {}_s\lambda_j^{(n)}\rceil/2, \lceil {}_s\lambda_j^{(n)}\rceil/2]\right) \cap \mathbf{Z}^d,$$

where $\lceil \cdot \rceil$ represents the "ceiling" function. Note that for $\mathbf{k} = (k_1, \ldots, k_d)' \in \mathbb{Z}^d \setminus E_n$, there exists $j \in \{1, \ldots, d\}$ such that $|k_j| > {}_s\lambda_j^{(n)}$, implying $\operatorname{dis}(R_{\mathbf{0},n} \cap \mathbf{Z}^d, R_{\mathbf{k},n} \cap \mathbf{Z}^d) \geq \operatorname{dis}(F_{\mathbf{0},n}, F_{\mathbf{k},n}) \geq 1$. Hence, sequentially using Lemmas 9.1 and 9.2, we may bound the covariances $\sigma_n(\mathbf{k})$, $\mathbf{k} \in \mathbb{Z}^d \setminus E_n$, with the mixing coefficient $\alpha(\cdot, \cdot)$,

$$|\sigma_n(\mathbf{k})| \leq 8[\mathrm{E}(Y_{\mathbf{0},n}^{2(2r+\delta)/r})]^{2r/(2r+\delta)} \alpha(\operatorname{dis}(R_{\mathbf{0},n} \cap \mathbf{Z}^d, R_{\mathbf{k},n} \cap \mathbf{Z}^d), {}_sN_n)^{\delta/(2r+\delta)}$$
$$\leq \mathcal{C}({}_sN_n)^{-2} \alpha(\operatorname{dis}(F_{\mathbf{0},n}, F_{\mathbf{k},n}), {}_sN_n)^{\delta/(2r+\delta)}.$$

From the above bound and $J_n(\mathbf{k})/|J_{\mathrm{OL}}| \leq 1$, $\mathbf{k} \in \mathbb{Z}^d$, we have

$$|W_{2n}| \leq \mathcal{C}|J_{\mathrm{OL}}|^{-1} \sum_{x=1}^{\infty} \left(\sum_{j=1}^{d} C_{(x,j,n)}\right) \alpha(x, {}_sN_n)^{\delta/(2r+\delta)},$$

(9.16)
$$C_{(x,j,n)} = |\{\mathbf{i} \in \mathbb{Z}^d : \operatorname{dis}(F_{\mathbf{0},n}, F_{\mathbf{i},n}) = x$$
$$= \inf\{|v_j - w_j| : \mathbf{v} \in F_{\mathbf{0},n}, \mathbf{w} \in F_{\mathbf{i},n}\}\}|.$$

The function $C_{(x,j,n)}$ counts the number of translated rectangles $F_{\mathbf{i},n}$ that lie a distance of $x \in \mathbb{Z}_+$ from the rectangle $F_{\mathbf{0},n}$, where this distance is realized in the $j$th coordinate direction for $j = 1, \ldots, d$. For $\mathbf{i} \in \mathbb{Z}^d$, $x \geq 1$ and $j \in \{1, \ldots, d\}$, if $\operatorname{dis}(F_{\mathbf{0},n}, F_{\mathbf{i},n}) = x = \inf\{|v_j - w_j| : \mathbf{v} \in F_{\mathbf{0},n}, \mathbf{w} \in F_{\mathbf{i},n}\}$, then $|i_j| = \lceil {}_s\lambda_j^{(n)} \rceil + x - 1$ with the remaining components of $\mathbf{i}$, namely $i_m$ for $m \in \{1, \ldots, d\} \setminus \{j\}$, constrained by $|i_m| \leq {}_s\lambda_m^{(n)} + x$. We use this observation to further bound the right-hand side of (9.16) by

$$\mathcal{C}|J_{\mathrm{OL}}|^{-1} \sum_{x=1}^{\infty} \left(\sum_{j=1}^{d} \prod_{m=1, j \neq m}^{d} 3({}_s\lambda_m^{(n)} + x)\right) \alpha(x, {}_sN_n)^{\delta/(2r+\delta)}$$
$$\leq \mathcal{C} \frac{\det({}_s\Delta_n)}{|J_{\mathrm{OL}}|} \sum_{j=1}^{d} ({}_s\lambda_{\min}^{(n)})^{-j} \left[\sum_{x=1}^{\ell_n} x^{j-1} + \sum_{x=\ell_n+1}^{\infty} x^{j-1} [\alpha_1(x) g({}_sN_n)]^{\delta/(2r+\delta)}\right]$$
$$\leq \mathcal{C} \frac{\det({}_s\Delta_n)}{|J_{\mathrm{OL}}|} \left[\frac{d\ell_n}{{}_s\lambda_{\min}^{(n)}} + \frac{\{\ell_n^{1/e}\}^{d\kappa\delta/(2r+\delta)}}{\ell_n^{2rd-d}} \sum_{x=\ell_n+1}^{\infty} x^{2rd-d-1} \alpha_1(x)^{\delta/(2r+\delta)}\right]$$
$$= o\left(\frac{\det({}_s\Delta_n)}{\det(\Delta_n)}\right),$$

using Assumptions A.1, A.3, Condition $M_r$ and $\ell_n = o({}_s\lambda_{\min}^{(n)})$ with $e$ from (9.1). This completes the proof of (9.13)(a).

To establish (9.13)(b), first note that

$$\left|\operatorname{Var}(\hat{\tau}_{n,\mathrm{OL}}^2) - \operatorname{Var}\left(\frac{{}_sN_n}{|J_{\mathrm{OL}}|} \sum_{\mathbf{i} \in J_{\mathrm{OL}}} Y_{\mathbf{i},n}^2\right)\right|$$

$$\leq 4\left(\sum_{j=1}^{5} A_{jn}^{1/2}\right)\left(\sum_{j=1}^{5} A_{jn}^{1/2} + \mathrm{Var}^{1/2}\left(\frac{{}_sN_n}{|J_{\mathrm{OL}}|}\sum_{\mathbf{i}\in J_{\mathrm{OL}}} Y_{\mathbf{i},n}^2\right)\right),$$

where $A_{1n} = \mathrm{Var}({}_sN_n \bar{Y}_n^2)$, $A_{2n} = \mathrm{Var}(|J_{\mathrm{OL}}|^{-1}{}_sN_n \sum_{\mathbf{i}\in J_{\mathrm{OL}}} Q_{\mathbf{i},n}^2)$, $A_{3n} = \mathrm{Var}({}_sN_n \bar{Q}_n^2)$, $A_{4n} = \mathrm{Var}(|J_{\mathrm{OL}}|^{-1}{}_sN_n \sum_{\mathbf{i}\in J_{\mathrm{OL}}} Y_{\mathbf{i},n}Q_{\mathbf{i},n})$, $A_{5n} = \mathrm{Var}({}_sN_n \bar{Y}_n \bar{Q}_n)$.

By (9.13)(a), it suffices to show that $A_{jn} = o(\det({}_s\Delta_n)/\det(\Delta_n))$ for each $j = 1, \ldots, 5$. We handle only two terms for illustration: $A_{1n}$, $A_{4n}$.

Consider $A_{1n}$. For $\mathbf{s} \in R_n \cap \mathbf{Z}^d$, let $\omega(\mathbf{s}) = [2^d \det({}_s\Delta_n)]^{-1}|\{\mathbf{i} \in J_{\mathrm{OL}} : \mathbf{s} \in \mathbf{i} + {}_s\Delta_n R_0\}|$ so that $0 \leq \omega(\mathbf{s}) \leq 1$. By Condition $M_r$ and Theorem 3 [Doukhan (1994), page 31] (similar to Lemma 9.2),

$$A_{1n} \leq \mathrm{E}(\bar{Y}_n^4)$$

(9.17)
$$= \frac{(2^d \det({}_s\Delta_n))^4}{|J_{\mathrm{OL}}|^4 ({}_sN_n)^2} \mathrm{E}\left(\left[\sum_{\mathbf{s}\in R_n \cap \mathbf{Z}^d} \omega(\mathbf{s})\nabla'(Z(\mathbf{s}) - \mu)\right]^4\right)$$

$$\leq \mathcal{C} \frac{(N_n)^2 (\det({}_s\Delta_n))^4}{|J_{\mathrm{OL}}|^4 ({}_sN_n)^2}.$$

Then $A_{1n} = o(\det({}_s\Delta_n)/\det(\Delta_n))$ follows from Lemma 9.3.

To handle $A_{4n}$, write $\sigma_{1n}(\mathbf{k}) = \mathrm{Cov}(Y_{\mathbf{0},n}Q_{\mathbf{0},n}, Y_{\mathbf{k},n}Q_{\mathbf{k},n})$, $\mathbf{k} \in \mathbb{Z}^d$. Then

$$A_{4n} = \frac{({}_sN_n)^2}{|J_{\mathrm{OL}}|^2} \sum_{\mathbf{k}\in\mathbb{Z}^d} J_n(\mathbf{k})\sigma_{1n}(\mathbf{k})$$

$$\leq \frac{({}_sN_n)^2}{|J_{\mathrm{OL}}|}\left(\sum_{\mathbf{k}\in E_n} |\sigma_{1n}(\mathbf{k})| + \sum_{\mathbf{k}\in\mathbb{Z}^d\setminus E_n} |\sigma_{1n}(\mathbf{k})|\right)$$

$$\equiv A_{4n}(E_n) + A_{4n}(E_n^c).$$

For $\mathbf{k} \in E_n$, note $|\sigma_{1n}(\mathbf{k})| \leq \mathcal{C}({}_sN_n)^{-3}$ using $|Y_{\mathbf{0},n}Q_{\mathbf{0},n}| \leq \mathcal{C}\|Z_{\mathbf{0},n} - \mu\|^3(1 + \|Z_{\mathbf{0},n} - \mu\|^a)$ (from Condition $D$) with Lemmas 9.1 and 9.2. From this bound, Lemma 9.3 and $|E_n| \leq 3^d \det({}_s\Delta_n)$, we find $A_{4n}(E_n) = o(\det({}_s\Delta_n)/\det(\Delta_n))$. We next bound the covariances $\sigma_{1n}(\mathbf{k})$, $\mathbf{k} \in \mathbb{Z}^d \setminus E_n$:

$$|\sigma_{1n}(\mathbf{k})| \leq 8[\mathrm{E}(|Y_{\mathbf{0},n}Q_{\mathbf{0},n}|^{(2r+\delta)/r})]^{2r/(2r+\delta)}$$

$$\times \alpha(\mathrm{dis}(R_{\mathbf{0},n} \cap \mathbf{Z}^d, R_{\mathbf{k},n} \cap \mathbf{Z}^d), {}_sN_n)^{\delta/(2r+\delta)}$$

$$\leq \mathcal{C}({}_sN_n)^{-3}\alpha(\mathrm{dis}(F_{\mathbf{0},n}, F_{\mathbf{k},n}), {}_sN_n)^{\delta/(2r+\delta)}$$

by the stationarity of the random field $Z(\cdot)$ and Lemmas 9.1 and 9.2. Using this inequality and repeating the same steps used to majorize "$W_{2n}$" from the proof of (9.13)(a) [see (9.16)], we have $A_{4n}(E_n^c) = o(\det({}_s\Delta_n)/\det(\Delta_n))$. The proof of Theorem 3.1(a) is complete.

9.2. *Proof of Theorem* 3.1(b). To simplify the counting arguments, we assume here integer-valued $_s\Delta_n \in \mathbb{Z}_+$, implying $_sN_{\mathbf{i},n} = {_sN_n}$, $\mathbf{i} \in \mathbb{Z}^d$. The more general case, in which the NOL subregions may differ in the number of sampling sites, is treated in Nordman (2002).

For each NOL subregion $\tilde{R}_{\mathbf{i},n}$, we denote the corresponding sample mean $\tilde{Z}_{\mathbf{i},n} = (_sN_{\mathbf{i},n})^{-1} \sum_{\mathbf{s} \in \tilde{R}_{\mathbf{i},n} \cap \mathbb{Z}^d} Z(\mathbf{s})$. The subsample evaluations of the statistic of interest, $\hat{\theta}_{\mathbf{i},n}^{\mathrm{NOL}}$, $\mathbf{i} \in J_{\mathrm{NOL}}$, can be expressed through a Taylor expansion of $H(\cdot)$ around $\mu$, substituting $\tilde{Z}_{\mathbf{i},n}$ for $Z_{\mathbf{i},n}$ in (9.12): $\hat{\theta}_{\mathbf{i},n}^{\mathrm{NOL}} = H(\tilde{Z}_{\mathbf{i},n}) = H(\mu) + \tilde{Y}_{\mathbf{i},n} + \tilde{Q}_{\mathbf{i},n}$.

We will complete the proof of Theorem 3.1(b) in two parts by showing

(9.18)
(a) $\quad \mathrm{Var}\left( \dfrac{_sN_n}{|J_{\mathrm{NOL}}|} \sum_{\mathbf{i} \in J_{\mathrm{NOL}}} \tilde{Y}_{\mathbf{i},n}^2 \right) = \dfrac{\det(_s\Delta_n)}{\det(\Delta_n)|R_0|} \cdot [2\tau^4](1 + o(1))$,

(b) $\quad \left| \mathrm{Var}(\hat{\tau}_{n,\mathrm{NOL}}^2) - \mathrm{Var}\left( \dfrac{_sN_n}{|J_{\mathrm{NOL}}|} \sum_{\mathbf{i} \in J_{\mathrm{NOL}}} \tilde{Y}_{\mathbf{i},n}^2 \right) \right| = o\left( \dfrac{\det(_s\Delta_n)}{\det(\Delta_n)} \right)$.

We will begin with showing (9.18)(a). For $\mathbf{k} \in \mathbb{Z}^d$, let $\tilde{J}_n(\mathbf{k}) = \{\mathbf{i} \in J_{\mathrm{NOL}} : \mathbf{i} + \mathbf{k} \in J_{\mathrm{NOL}}\}$ and $\tilde{\sigma}_n(\mathbf{k}) = \mathrm{Cov}(\tilde{Y}_{\mathbf{0},n}^2, \tilde{Y}_{\mathbf{k},n}^2)$. Then we may express the variance,

(9.19)
$$\mathrm{Var}\left( \dfrac{_sN_n}{|J_{\mathrm{NOL}}|} \sum_{\mathbf{i} \in J_{\mathrm{NOL}}} \tilde{Y}_{\mathbf{i},n}^2 \right)$$
$$= \dfrac{(_sN_n)^2}{|J_{\mathrm{NOL}}|^2} \left( \sum_{\mathbf{k} \in \mathbb{Z}^d, 0 < \|\mathbf{k}\|_\infty \leq 1} \tilde{J}_n(\mathbf{k}) \tilde{\sigma}_n(\mathbf{k}) \right.$$
$$\left. + \sum_{\mathbf{k} \in \mathbb{Z}^d, \|\mathbf{k}\|_\infty > 1} \tilde{J}_n(\mathbf{k}) \tilde{\sigma}_n(\mathbf{k}) + |J_{\mathrm{NOL}}| \tilde{\sigma}_n(\mathbf{0}) \right)$$
$$\equiv U_{1n} + U_{2n} + |J_{\mathrm{NOL}}|^{-1} (_sN_n)^2 \tilde{\sigma}_n(\mathbf{0}).$$

We first prove $U_{2n} = o(|J_{\mathrm{NOL}}|^{-1})$, noting that $\det(_s\Delta_n)/\det(\Delta_n) = O(|J_{\mathrm{NOL}}|^{-1})$ from Lemma 9.3.

When $\mathbf{k} = (k_1, \ldots, k_d)' \in \mathbb{Z}^d$, $\|\mathbf{k}\|_\infty > 1$, then for some $1 \leq m_k \leq d$,

$$\mathrm{dis}(\tilde{R}_{\mathbf{0},n} \cap \mathbf{Z}^d, \tilde{R}_{\mathbf{k},n} \cap \mathbf{Z}^d) \geq \max_{1 \leq j \leq d}(|k_j| - 1) {_s\lambda_j^{(n)}}$$
$$\equiv (|k_{m_k}| - 1) {_s\lambda_{m_k}^{(n)}}.$$

If $j \in \{1, \ldots, d\}, j \neq m_k$, we have

$$|k_j| \leq (_s\lambda_j^{(n)})^{-1} (|k_{m_k}| - 1) {_s\lambda_{m_k}^{(n)}} + 1.$$

Note also if $\mathbf{k} \in \mathbb{Z}^d$, $\|\mathbf{k}\|_\infty > 1$, then

$$|\tilde{\sigma}(\mathbf{k})| \leq \mathcal{C}(_sN_n)^{-2}\alpha((|k_{m_k}|-1)_s\lambda_{m_k}^{(n)}, {}_sN_n)$$

by Lemmas 9.1 and 9.2. Hence, we have

$$|U_{2n}| \leq \frac{\mathcal{C}}{|J_{\text{NOL}}|} \sum_{x=1}^{\infty} \left[ \sum_{j=1}^{d} \prod_{i=1, i\neq j}^{d} (_s\lambda_i^{(n)} x + {}_s\lambda_j^{(n)})/{}_s\lambda_i^{(n)} \right]$$

$$\times \alpha(_s\lambda_{\min}^{(n)} x, {}_sN_n)^{\delta/(2r+\delta)}$$

$$\leq \frac{\mathcal{C}}{|J_{\text{NOL}}|} \sum_{x=1}^{\infty} x^{d-1} \left( \frac{{}_s\lambda_n^{\max}}{{}_s\lambda_{\min}^{(n)}} \right)^{d-1} \alpha(_s\lambda_{\min}^{(n)} x, {}_sN_n)^{\delta/(2r+\delta)}$$

$$\leq \frac{\mathcal{C}}{|J_{\text{NOL}}|} \frac{\{\det(_s\Delta_n)\}^{\kappa\delta/(2r+\delta)}}{(_s\lambda_{\min}^{(n)})^{2rd-d-1}} \sum_{x=1}^{\infty} (_s\lambda_{\min}^{(n)} x)^{2rd-d-1} \alpha_1(_s\lambda_{\min}^{(n)} x)^{\delta/(2r+\delta)}$$

$$= o\left(\frac{1}{|J_{\text{NOL}}|}\right)$$

by Assumptions A.1, A.3 and Condition $M_r$.

We now show that $U_{1n} = o(|J_{\text{NOL}}|^{-1})$. For $\mathbf{k} \in \mathbb{Z}^d$, $0 < \|\mathbf{k}\|_\infty \leq 1$, define the set

$$T_{\mathbf{k},n}^j =: \begin{cases} \{\mathbf{x} \in \mathbb{R}^d : 1/2 \cdot {}_s\lambda_j^{(n)} < x_j \leq 1/2 \cdot {}_s\lambda_j^{(n)} + \ell_n\}, & \text{if } k_j = 1, \\ \{\mathbf{x} \in \mathbb{R}^d : -1/2 \cdot {}_s\lambda_j^{(n)} - \ell_n < x_j \leq -1/2 \cdot {}_s\lambda_j^{(n)}\}, & \text{if } k_j = -1, \\ \varnothing, & \text{if } k_j = 0, \end{cases}$$

for each coordinate direction $j = 1, \ldots, d$. Let $T_{\mathbf{k},n} = \bigcup_{j=1}^{d} T_{\mathbf{k},n}^j$. We decompose the sum: ${}_sN_n \tilde{Y}_{\mathbf{k},n} = \overline{\Sigma}(\tilde{R}_{\mathbf{k},n} \setminus T_{\mathbf{k},n}) + \overline{\Sigma}(\tilde{R}_{\mathbf{k},n} \cap T_{\mathbf{k},n}) \equiv \tilde{S}_{\mathbf{k},n} + \tilde{S}_{\mathbf{k},n}^*$. Then, $U_{1n} = o(|J_{\text{NOL}}|^{-1})$ follows from 1–4 below:

1. $|\mathrm{E}(\tilde{Y}_{\mathbf{0},n}^2 \tilde{S}_{\mathbf{k},n} \tilde{S}_{\mathbf{k},n}^*)| \leq [\mathrm{E}(\tilde{Y}_{\mathbf{0},n}^6) \mathrm{E}(|\tilde{S}_{\mathbf{k},n}|^3) \mathrm{E}(|\tilde{S}_{\mathbf{k},n}^*|^3)]^{1/3} = o(1)$, using Lemma 9.2 and

$$|\tilde{R}_{\mathbf{k},n} \cap T_{\mathbf{k},n} \cap \mathbb{Z}^d| \leq \sum_{j=1}^{d} |\tilde{R}_{\mathbf{k},n} \cap T_{\mathbf{k},n}^j \cap \mathbb{Z}^d| \leq \ell_n \det(_s\Delta_n) \sum_{j=1}^{d} (_s\lambda_j^{(n)})^{-1}$$

$$= o(_s\lambda_{\min}^{(n)}).$$

2. Likewise, $\mathrm{E}(\tilde{Y}_{\mathbf{0},n}^2 \tilde{S}_{\mathbf{k},n}^{*2}) \leq [\mathrm{E}(\tilde{Y}_{\mathbf{0},n}^4) \mathrm{E}(\tilde{S}_{\mathbf{k},n}^{*4})]^{1/2} = o(1)$.
3. $|{}_sN_n \mathrm{E}(\tilde{Y}_{\mathbf{k},n}^2) - ({}_sN_n)^{-1} \mathrm{E}(\tilde{S}_{\mathbf{k},n}^2)| \leq 4({}_sN_n)^{-1} \max\{[\mathrm{E}(\tilde{S}_{\mathbf{k},n}^2)\mathrm{E}(\tilde{S}_{\mathbf{k},n}^{*2})]^{1/2}, \mathrm{E}(\tilde{S}_{\mathbf{k},n}^{*2})\} = o(1)$.
4. $|\mathrm{Cov}(\tilde{Y}_{\mathbf{0},n}^2, \tilde{S}_{\mathbf{k},n}^2)| \leq \mathcal{C}\alpha(\ell_n, {}_sN_n)^{\delta/(2r+\delta)} = o(1)$ by applying Lemmas 9.1 and 9.2, Assumption A.3, and Condition $M_r$ with $\mathrm{dis}(\tilde{R}_{\mathbf{k},n} \cap \mathbb{Z}^d \setminus T_{\mathbf{k},n}, {}_sR_n \cap \mathbb{Z}^d) \geq \ell_n$.

Since $\tilde{\sigma}_n(\mathbf{0}) = \text{Var}(Y_{\mathbf{0},n}^2)$, the remaining quantity in (9.19) can be expressed as

$$\frac{({}_sN_n)^2}{|J_{\text{NOL}}|}\tilde{\sigma}_n(\mathbf{0}) = \frac{1}{|J_{\text{NOL}}|}\text{Var}([\nabla' Z_\infty]^2)(1+o(1))$$

$$= \frac{\det({}_s\Delta_n)}{\det(\Delta_n)|R_0|} \cdot [2\tau^4](1+o(1))$$

by applying the CLT [as in (9.9)] and Lemma 9.3. We have now established (9.18)(a).

We omit the proof of (9.18)(b), which resembles the one establishing (9.13)(b) and incorporates arguments used to bound $U_{1n}, U_{2n}$; Nordman (2002) provides more details.

**10. Proofs for bias expansions.** We will use the following lemma concerning $\tau_n^2 = N_n \text{Var}(\hat{\theta}_n)$ to prove the theorems pertaining to bias expansions of $\hat{\tau}_{n,\text{OL}}^2$ and $\hat{\tau}_{n,\text{NOL}}^2$.

LEMMA 10.1. *Under the assumptions and conditions of Theorem* 3.1,

$$\tau_n^2 = \tau^2 + O([\det(\Delta_n)]^{-1/\max\{2,d\}}).$$

PROOF. By a Taylor expansion around $\mu$: $\hat{\theta}_n = H(\bar{Z}_{N_n}) = H(\mu) + \bar{Y}_{N_n} + \bar{Q}_{N_n}$ [replacing $\bar{Z}_{N_n}$ for $Z_{\mathbf{i},n}$ in (9.12)] and so $N_n\text{Var}(\hat{\theta}_n) = N_n\text{Var}(\bar{Y}_{N_n} + \bar{Q}_{N_n})$. For $\mathbf{k} \in \mathbb{Z}^d$, let $N_n(\mathbf{k}) = |\{\mathbf{i} \in R_n \cap \mathbf{Z}^d : \mathbf{i} + \mathbf{k} \in R_n\}|$. It holds that $N_n(\mathbf{k}) \leq N_n$ and

$$N_n \leq N_n(\mathbf{k}) + |\{\mathbf{i} \in \mathbb{Z}^d : T^{\mathbf{i}} \cap \overline{R_n} \neq \varnothing, \ T^{\mathbf{i}} \cap \overline{R_n^c} \neq \varnothing;$$
$$(10.1) \qquad\qquad T^{\mathbf{i}} = \mathbf{i} + \|\mathbf{k}\|_\infty [-1,1]^d\}|$$
$$\leq N_n(\mathbf{k}) + \mathcal{C}\|\mathbf{k}\|_\infty^d (\lambda_n^{\max})^{d-1}$$

by the boundary condition on $R_0$. Also, by Lemma 9.1 and stationarity, for each $\mathbf{k} \neq \mathbf{0} \in \mathbb{Z}^d$,

$$(10.2) \qquad |\sigma(\mathbf{k})| \leq \mathcal{C}\alpha_1(\|\mathbf{k}\|_\infty)^{\delta/(2r+\delta)}, \qquad \mathbf{k} \in \mathbb{Z}^d.$$

Using $|\{\mathbf{k} \in \mathbb{Z}^d : \|\mathbf{k}\|_\infty = x\}| \leq \mathcal{C}x^{d-1}$, $x \geq 1$, the covariances are absolutely summable over $\mathbb{Z}^d$:

$$(10.3) \qquad \sum_{\mathbf{k} \in \mathbb{Z}^d} |\sigma(\mathbf{k})| \leq |\sigma(\mathbf{0})| + \mathcal{C}\sum_{x=1}^\infty x^{d-1}\alpha_1(x)^{\delta/(2r+\delta)} < \infty.$$

From (10.1)–(10.3), we find

$$(10.4) \quad N_n \operatorname{Var}(\bar{Y}_{N_n}) = \frac{1}{N_n} \sum_{\mathbf{k} \in \mathbb{Z}^d} N_n(\mathbf{k}) \sigma(\mathbf{k}) = \tau^2 + I_n,$$

$$|I_n| \leq \frac{1}{N_n} \sum_{\mathbf{k} \in \mathbb{Z}^d} |N_n - N_n(\mathbf{k})| \cdot |\sigma(\mathbf{k})|$$

$$(10.5) \quad \leq \mathcal{C} \cdot \frac{(\lambda_n^{\max})^{d-1}}{N_n} \sum_{x=1}^{\infty} x^{2d-1} \alpha_1(x)^{\delta/(2r+\delta)}$$

$$= O([\det(\Delta_n)]^{-1/d}).$$

By Condition $D$ and Lemma 9.2, it follows that $N_n \operatorname{Var}(\bar{Q}_{N_n}) = O([\det(\Delta_n)]^{-1})$.

Finally, with bounds on the variance of $\bar{Y}_{N_n}$ and $\bar{Q}_{N_n}$, we apply the Cauchy–Schwarz inequality to the covariance $N_n |\operatorname{Cov}(\bar{Y}_{N_n}, \bar{Q}_{N_n})| = O([\det(\Delta_n)]^{-1/2})$, setting the order on the difference $|N_n \operatorname{Var}(\hat{\theta}_n) - \tau^2|$. □

We give a few lemmas which help compute the bias of the estimators $\hat{\tau}^2_{n,\mathrm{OL}}$ and $\hat{\tau}^2_{n,\mathrm{NOL}}$.

LEMMA 10.2. *Let $\tilde{Y}_{\mathbf{i},n} = ({}_sN_{\mathbf{i},n})^{-1} \sum_{\mathbf{s} \in \mathbf{Z}^d \cap \tilde{R}_{\mathbf{i},n}} \nabla'(Z(\mathbf{s}) - \mu)$, $\mathbf{i} \in \mathbb{Z}^d$. Suppose Assumptions A.1–A.5 and Conditions $D_2$ and $M_{2+a}$ hold with $d \geq 2$ with $a$ as specified under Condition $D_2$. Then*

$$\mathrm{E}(\hat{\tau}^2_{n,\mathrm{OL}}) - {}_sN_{\mathbf{0},n} \mathrm{E}(\tilde{Y}^2_{\mathbf{0},n}), \mathrm{E}(\hat{\tau}^2_{n,\mathrm{NOL}}) - |J_{\mathrm{NOL}}|^{-1} \sum_{\mathbf{i} \in J_{\mathrm{NOL}}} {}_sN_{\mathbf{i},n} \mathrm{E}(\tilde{Y}^2_{\mathbf{i},n})$$

$$=: O([\det({}_s\Delta_n)]^{-1/2}) + o([\det({}_s\Delta_n)]^{-1/d}).$$

PROOF. We consider here only $\mathrm{E}(\hat{\tau}^2_{n,\mathrm{OL}})$. For integer ${}_s\Delta_n$, the arguments for $\mathrm{E}(\hat{\tau}^2_{n,\mathrm{NOL}})$ are essentially the same; more details are provided in Nordman (2002).

By stationarity and an algebraic expansion as in (9.12),

$$\mathrm{E}(\hat{\tau}^2_{n,\mathrm{OL}}) = {}_sN_n [\mathrm{E}(Y^2_{\mathbf{0},n}) + \mathrm{E}(Q^2_{\mathbf{0},n})$$
$$+ 2\mathrm{E}(Y_{\mathbf{0},n} Q_{\mathbf{0},n}) - \mathrm{E}(\bar{Y}^2_n) - \mathrm{E}(\bar{Q}^2_n) - 2\mathrm{E}(\bar{Y}_n \bar{Q}_n)].$$

With the moment arguments based on Lemma 9.2 and Condition $D_r$, we have

$$(10.6) \quad \begin{aligned} {}_sN_n \mathrm{E}(Y^2_{\mathbf{0},n}) &\leq \mathcal{C}, \\ {}_sN_n \mathrm{E}(\bar{Y}^2_n) &\leq \mathcal{C} {}_sN_n(N_n)^{-1}, \\ {}_sN_n \mathrm{E}(Q^2_{\mathbf{0},n}), \quad {}_sN_n \mathrm{E}(\bar{Q}^2_n) &\leq \mathcal{C}({}_sN_n)^{-1}, \end{aligned}$$

where the bound on ${}_sN_n\mathrm{E}(\bar{Y}_n^2)$ follows from (9.17). By Hölder's inequality and Assumption A.2,

$$\mathrm{E}(\hat{\tau}_{n,\mathrm{OL}}^2) = {}_sN_n\mathrm{E}(Y_{\mathbf{0},n}^2) + O(({}_sN_n)^{-1/2}) + O({}_sN_n(N_n)^{-1}).$$

Note that $\mathrm{E}(Y_{\mathbf{0},n}^2) = \mathrm{E}(\tilde{Y}_{\mathbf{0},n}^2)$, ${}_sN_n = {}_sN_{\mathbf{0},n}$. Hence, applying Lemma 9.3 and Assumption A.2, we establish Lemma 10.2 for $\hat{\tau}_{n,\mathrm{OL}}^2$. □

The next lemma provides a small refinement to Lemma 10.2 made possible when the function $H(\cdot)$ is smoother. We shall make use of this lemma in bias expansions of $\hat{\tau}_{n,\mathrm{OL}}^2$ and $\hat{\tau}_{n,\mathrm{NOL}}^2$ in lower sampling dimensions, namely $d = 1$ or 2.

LEMMA 10.3. *Assume $d = 1$ or 2. In addition to Assumptions A.1–A.5, suppose that Conditions $D_3$ and $M_{3+a}$ hold with $a$ as specified under Condition $D_3$. Then*

$$\mathrm{E}(\hat{\tau}_{n,\mathrm{OL}}^2) - {}_sN_{\mathbf{0},n}\mathrm{E}(\tilde{Y}_{\mathbf{0},n}^2), \mathrm{E}(\hat{\tau}_{n,\mathrm{NOL}}^2) - |J_{\mathrm{NOL}}|^{-1}\sum_{\mathbf{i}\in J_{\mathrm{NOL}}} {}_sN_{\mathbf{i},n}\mathrm{E}(\tilde{Y}_{\mathbf{i},n}^2)$$

$$=: \begin{cases} O([\det({}_s\Delta_n)]^{-1}), & \text{if } d = 1, \\ o([\det({}_s\Delta_n)]^{-1/2}), & \text{if } d = 2. \end{cases}$$

PROOF. We again consider only $\hat{\tau}_{n,\mathrm{OL}}^2$. For $\mathbf{i} \in J_{\mathrm{OL}}$, we use a third-order Taylor expansion of each subsample statistic around $\mu$: $\hat{\theta}_{\mathbf{i},n} = H(\mu) + Y_{\mathbf{i},n} + Q_{\mathbf{i},n} + C_{\mathbf{i},n}$, where $Y_{\mathbf{i},n} = \nabla'(Z_{\mathbf{i},n} - \mu)$,

$$Q_{\mathbf{i},n} = \sum_{\|\alpha\|_1 = 2} \frac{c_\alpha}{\alpha!}(Z_{\mathbf{i},n} - \mu)^\alpha,$$

$$C_{\mathbf{i},n} = 3\sum_{\|\alpha\|_1 = 3} \frac{c_\alpha}{\alpha!}(Z_{\mathbf{i},n} - \mu)^\alpha \int_0^1 (1-\omega)^2 D^\alpha H(\mu + \omega(Z_{\mathbf{i},n} - \mu))\,d\omega.$$

Here $C_{\mathbf{i},n}$ denotes the remainder term in the Taylor expansion and $Q_{\mathbf{i},n}$ is defined a little differently here compared to (9.12). Write the sample means for the Taylor terms: $\bar{Y}_n$, $\bar{Q}_n$ as before, $\bar{C}_n = |J_{\mathrm{OL}}|^{-1}\sum_{\mathbf{i}\in J_{\mathrm{OL}}} C_{\mathbf{i},n}$. The moment inequalities in (10.6) are still valid and, by Lemma 9.2 and Condition $D$, we can produce bounds ${}_sN_n\mathrm{E}(C_{\mathbf{0},n}^2), {}_sN_n\mathrm{E}(\bar{C}_n^2) \leq \mathcal{C}({}_sN_n)^{-2}$. By Hölder's inequality and the scaling conditions from Assumptions A.1 and A.2, we then have

$$\mathrm{E}(\hat{\tau}_{n,\mathrm{OL}}^2) = {}_sN_n[\mathrm{E}(Y_{\mathbf{0},n}^2) + 2\mathrm{E}(Y_{\mathbf{0},n}Q_{\mathbf{0},n})] + \begin{cases} O([\det({}_s\Delta_n)]^{-1}), & \text{if } d = 1, \\ o([\det({}_s\Delta_n)]^{-1/2}), & \text{if } d = 2. \end{cases}$$

Since ${}_sN_n\mathrm{E}(Y_{\mathbf{0},n}^2) = {}_sN_{\mathbf{0},n}\mathrm{E}(\tilde{Y}_{\mathbf{0},n}^2)$, Lemma 10.3 for $\hat{\tau}_{n,\mathrm{OL}}^2$ will follow by showing

$${}_sN_n\mathrm{E}(Y_{\mathbf{0},n}Q_{\mathbf{0},n})$$
$$(10.7) \qquad = {}_sN_n \sum_{i,j,k=1}^{p} c_i a_{j,k} \mathrm{E}[(Z_{i,\mathbf{0},n} - \mu)(Z_{j,\mathbf{0},n} - \mu)(Z_{k,\mathbf{0},n} - \mu)]$$
$$= O([\det({}_s\Delta_n)]^{-1}),$$

where $Z_{\mathbf{0},n} = (Z_{1,\mathbf{0},n}, \ldots, Z_{p,\mathbf{0},n})' \in \mathbb{R}^p$ is a vector of coordinate sample means, $c_i = \partial H(\mu)/\partial x_i$; $a_{j,k} = 1/2 \cdot \partial^2 H(\mu)/\partial x_j \partial x_k$.

Denote the observation $Z(\mathbf{s}) = (Z_1(\mathbf{s}), \ldots, Z_p(\mathbf{s}))' \in \mathbb{R}^d$, $\mathbf{s} \in \mathbf{Z}^d$. Fix $i,j,k \in \{1,\ldots,p\}$ and w.l.o.g. assume $\mu = 0$. Then ${}_sN_n|\mathrm{E}(Z_{i,\mathbf{0},n}Z_{j,\mathbf{0},n}Z_{k,\mathbf{0},n})| = |({}_sN_n)^{-1}\mathrm{E}(Z_i(\mathbf{t})Z_j(\mathbf{t})Z_k(\mathbf{t})) + L_{1n}^{ijk} + L_{2n}^{ijk}|$ where

$$L_{1n}^{ijk} = ({}_sN_n)^{-2} \sum_{\substack{\mathbf{u},\mathbf{v},\mathbf{w}\in\mathbf{Z}^d\cap{}_sR_n \\ \mathbf{u}\neq\mathbf{v}\neq\mathbf{w}}} \mathrm{E}[Z_i(\mathbf{u})Z_j(\mathbf{v})Z_k(\mathbf{w})],$$

$$L_{2n}^{ijk} = ({}_sN_n)^{-2} \sum_{\substack{\mathbf{u},\mathbf{v}\in\mathbf{Z}^d\cap{}_sR_n \\ \mathbf{u}\neq\mathbf{v}}} \mathrm{E}[Z_i(\mathbf{u})Z_j(\mathbf{u})Z_k(\mathbf{v})$$
$$+ Z_i(\mathbf{u})Z_j(\mathbf{v})Z_k(\mathbf{u}) + Z_i(\mathbf{v})Z_j(\mathbf{u})Z_k(\mathbf{u})].$$

By Lemma 9.1, Assumption A.3 and Condition $M_r$,

$$|L_{2n}^{ijk}| \leq \frac{\mathcal{C}}{{}_sN_n} \sum_{x=1}^{\infty} x^{d-1}\alpha(x,1)^{\delta/(2r+\delta)} = O([\det({}_s\Delta_n)]^{-1}),$$

similarly to (10.3). For $\mathbf{y}_1, \mathbf{y}_2, \mathbf{y}_3 \in \mathbb{R}^d$, define $\mathrm{dis}_3(\{\mathbf{y}_1,\mathbf{y}_2,\mathbf{y}_3\}) = \max_{1\leq i\leq 3}\mathrm{dis}(\{\mathbf{y}_i\},\{\mathbf{y}_1,\mathbf{y}_2,\mathbf{y}_3\}\setminus\{\mathbf{y}_i\})$. If $x \geq 1 \in \mathbb{Z}_+$, then $|\{(\mathbf{y}_1,\mathbf{y}_2) \in (\mathbb{Z}^d)^2 : \mathrm{dis}_3(\{\mathbf{y}_1,\mathbf{y}_2,\mathbf{0}\}) = x\}| \leq \mathcal{C}x^{2d-1}$ from Theorem 4.1, Lahiri (1999a). Thus,

$$|L_{1n}^{ijk}| \leq \frac{\mathcal{C}}{{}_sN_n} \sum_{x=1}^{\infty} x^{2d-1}\alpha(x,2)^{\delta/(2r+\delta)} = O([\det({}_s\Delta_n)]^{-1}).$$

This establishes (10.7), completing the proof of Lemma 10.3 for $\hat{\tau}_{n,\mathrm{OL}}^2$. $\square$

We use the next lemma in the proof of Theorem 4.2. It allows us to approximate lattice point counts with Lebesgue volumes, in $\mathbb{R}^2$ or $\mathbb{R}^3$, to a sufficient degree of accuracy.

LEMMA 10.4. *Let $d = 2, 3$ and $R_0 \subset (-1/2, 1/2]^d$ such that $B^\circ \subset R_0 \subset \overline{B}$ for a convex set $B$. Let $\{b_n\}_{n=1}^{\infty}$ be a sequence of positive real numbers such* that $b_n \to \infty$. If $\mathbf{k} \in \mathbb{Z}^d$, then there exist $N_\mathbf{k} \in \mathbb{Z}_+$ and $\mathcal{C}_d > 0$ such that for $n \geq N_\mathbf{k}$, $\mathbf{i} \in \mathbb{Z}^d$,

$$|(|b_n R_0| - |\mathbf{Z}^d \cap b_n(\mathbf{i} + R_0)|)$$
$$- (|b_n R_0 \cap \mathbf{k} + b_n R_0| - |\mathbf{Z}^d \cap b_n(\mathbf{i} + R_0) \cap \mathbf{k} + b_n(\mathbf{i} + R_0)|)|$$
$$\leq \begin{cases} \mathcal{C}_2 \|\mathbf{k}\|_\infty^2, & \text{if } d = 2, \\ \mathcal{C}_3 \|\mathbf{k}\|_\infty^4 (b_n^{5/3} + \xi_{\mathbf{k},n} b_n^2), & \text{if } d = 3, \end{cases}$$

where $\{\xi_{\mathbf{k},n}\}_{n=1}^\infty \subset \mathbb{R}$ is a nonnegative sequence (possibly dependent on $\mathbf{k}$) such that $\xi_{\mathbf{k},n} \to 0$.

The proof is provided in Nordman and Lahiri (2002).

To establish Lemma 4.1, we require some additional notation. For $\mathbf{i}, \mathbf{k} \in \mathbb{Z}^d$, and let $_sN_{\mathbf{i},n}(\mathbf{k}) = |\mathbf{Z}^d \cap \tilde{R}_{\mathbf{i},n} \cap \mathbf{k} + \tilde{R}_{\mathbf{i},n}|$ denote the number of sampling sites (lattice points) in the intersection of a NOL subregion with its $\mathbf{k}$-translate. Note $_sN_{\mathbf{i},n}(\mathbf{k})$ is a subsample version of $N_n(\mathbf{k})$ from (10.1).

PROOF OF LEMMA 4.1. We start with bounds

(10.8) $$\sup_{\mathbf{i} \in \mathbb{Z}^d} |_sN_n - {}_sN_{\mathbf{i},n}| \leq \mathcal{C}({}_s\lambda_n^{\max})^{d-1},$$

(10.9) $$\begin{aligned} |_sN_{\mathbf{i},n} - {}_sN_{\mathbf{i},n}(\mathbf{k})| \\ \leq |\{\mathbf{j} \in \mathbb{Z}^d : T^{\mathbf{j}} \cap \overline{{}_sR_n} \neq \varnothing,\ T^{\mathbf{j}} \cap \overline{{}_sR_n^c} \neq \varnothing;\ T^{\mathbf{j}} = \mathbf{j} + \|\mathbf{k}\|_\infty [-2, 2]^d\}| \\ \leq \mathcal{C}\|\mathbf{k}\|_\infty^d ({}_s\lambda_n^{\max})^{d-1}, \end{aligned}$$

by the boundary condition on $R_0$ (cf. Lemma 9.3) and $\inf_{\mathbf{j} \in \mathbb{Z}^d} \|{}_s\Delta_n \mathbf{i} - \mathbf{j}\|_\infty \leq 1/2$.

Modify (10.4) by replacing $N_n, N_n(\mathbf{k}), \bar{Y}_{N_n}$ with $_sN_{\mathbf{i},n}, {}_sN_{\mathbf{i},n}(\mathbf{k}), \tilde{Y}_{\mathbf{i},n} = \nabla'(\tilde{Z}_{\mathbf{i},n} - \mu)$ (i.e., use a NOL subregion in place of the sampling region), and replace $N_n, \Delta_n, \lambda_n^{\max}$ with the subsample analogs $_sN_{\mathbf{i},n}, {}_s\Delta_n, {}_s\lambda_n^{\max}$ in (10.5). We then find, using (10.3), for each $\mathbf{i} \in \mathbb{Z}^d$,

(10.10) $$_sN_{\mathbf{i},n} \mathrm{E}(\tilde{Y}_{\mathbf{i},n}^2) - \tau^2 = \frac{1}{{}_sN_{\mathbf{i},n}} \sum_{\mathbf{k} \in \mathbb{Z}^d} ({}_sN_{\mathbf{i},n}(\mathbf{k}) - {}_sN_{\mathbf{i},n})\sigma(\mathbf{k}) \equiv {}_sI_{\mathbf{i},n},$$

(10.11) $$\begin{aligned} \sup_{\mathbf{i} \in \mathbb{Z}^d} |_sI_{\mathbf{i},n}| &\leq \sup_{\mathbf{i} \in \mathbb{Z}^d} \left\{ \frac{1}{{}_sN_{\mathbf{i},n}} \sum_{\mathbf{k} \in \mathbb{Z}^d} |_sN_{\mathbf{i},n}(\mathbf{k}) - {}_sN_{\mathbf{i},n}| \cdot |\sigma(\mathbf{k})| \right\} \\ &\leq \mathcal{C} \cdot \frac{({}_s\lambda_n^{\max})^{d-1}}{{}_sN_n - \mathcal{C}({}_s\lambda_n^{\max})^{d-1}} \sum_{x=1}^\infty x^{2d-1} \alpha_1(x)^{\delta/(2r+\delta)} \\ &= O([\det({}_s\Delta_n)]^{-1/d}), \end{aligned}$$

OPTIMAL SPATIAL SUBSAMPLE SIZE 45from (10.8), (10.9) and Assumption A.1. Now applying Lemma 10.1 and Assumption A.2 with Lemma 10.2 for $d \geq 2$ or Lemma 10.3 for $d = 1$, Lemma 4.1 follows. $\square$

PROOF OF THEOREM 4.1. Here $_sN_{\mathbf{i},n} = {_sN_n}$, $_sN_{\mathbf{i},n}(\mathbf{k}) = C_n(\mathbf{k})$, $\mathrm{E}(\tilde{Y}_{\mathbf{i},n}) = \mathrm{E}(\tilde{Y}_{\mathbf{0},n})$ for each $\mathbf{i}, \mathbf{k} \in \mathbb{Z}^d$ (since $_s\lambda_n \in \mathbb{Z}_+$ for NOL subsamples) and $\det(_s\Delta_n) = {_s\lambda_n}^d$. Applying Lemma 10.2 for $d \geq 3$ and Lemma 10.3 for $d = 2$, Lemma 10.1, Assumption A.2 and (10.11),

$$\mathrm{E}(\hat{\tau}_n^2) - \tau_n^2 = \frac{-1}{{_s\lambda_n}|R_0|} \sum_{\mathbf{k} \in \mathbb{Z}^d} g_n(\mathbf{k}) + o({_s\lambda_n}^{-1}),$$

$$g_n(k) \equiv \frac{{_sN_n} - C_n(\mathbf{k})}{{_s\lambda_n}^{d-1}} \cdot \frac{{_s\lambda_n}^d |R_0|}{{_sN_n}} \cdot \sigma(\mathbf{k}).$$

From (10.11) and Lemma 9.3, it follows that $\sum_{\mathbf{k} \in \mathbb{Z}^d} |g_n(\mathbf{k})| \leq \mathcal{C}$, $n \in \mathbb{Z}_+$, and that $g_n(\mathbf{k}) \to C(\mathbf{k})\sigma(\mathbf{k})$ for $\mathbf{k} \in \mathbb{Z}^d$. By the LDCT, the proof of Theorem 4.1 is complete. $\square$

To establish Theorem 4.2, we require some additional notation. For $\mathbf{i}, \mathbf{k} \in \mathbb{Z}^d$, denote the difference between two Lebesgue volume-for-count approximations as

$$D_{\mathbf{i},n}(\mathbf{k}) = (|\tilde{R}_{\mathbf{i},n}| - {_sN_{\mathbf{i},n}}) - (|\tilde{R}_{\mathbf{i},n} \cap \mathbf{k} + \tilde{R}_{\mathbf{i},n}| - {_sN_{\mathbf{i},n}}(\mathbf{k}))$$
$$= (|_sR_n| - {_sN_{\mathbf{i},n}}) - (|_sR_n \cap \mathbf{k} + {_sR_n}| - {_sN_{\mathbf{i},n}}(\mathbf{k})).$$

PROOF OF THEOREM 4.2. We handle here the cases $d = 2$ or $3$. Details on the proof for $d = 1$ are given in Nordman (2002). We note first that if $V(\mathbf{k})$ exists for each $\mathbf{k} \in \mathbb{Z}^d$ then Lemma 10.4 implies $C(\mathbf{k}) = V(\mathbf{k})$.

Consider $\hat{\tau}_{n,\mathrm{NOL}}^2$. Applying Lemma 10.2 for $d = 3$, and Lemma 10.3 for $d = 2$, with (10.8), (10.10) and (10.11) gives

$$\mathrm{E}(\hat{\tau}_{n,\mathrm{NOL}}^2) - \tau_n^2 = |J_{\mathrm{NOL}}|^{-1} \sum_{\mathbf{i} \in J_{\mathrm{NOL}}} \frac{{_sN_{\mathbf{i},n}}}{|_sR_n|} {_sI_{\mathbf{i},n}} + o({_s\lambda_n}^{-1}).$$

Then, using (10.3), we can arrange terms to write

$$|J_{\mathrm{NOL}}|^{-1} \sum_{\mathbf{i} \in J_{\mathrm{NOL}}} \frac{{_sN_{\mathbf{i},n}}}{|_sR_n|} {_sI_{\mathbf{i},n}} = \Psi_n + \sum_{\mathbf{k} \in \mathbb{Z}^d} \frac{G_n(\mathbf{k})}{{_s\lambda_n}|R_0|};$$

$$G_n(\mathbf{k}) = \sum_{\mathbf{i} \in J_{\mathrm{NOL}}} \frac{D_{\mathbf{i},n}(\mathbf{k})\sigma(\mathbf{k})}{{_s\lambda_n}^{d-1}|J_{\mathrm{NOL}}|}$$

for $\Psi_n = -\sum_{\mathbf{k} \in \mathbb{Z}^d} |{}_sR_n|^{-1}(|{}_sR_n| - |{}_sR_n \cap \mathbf{k} + {}_sR_n|)\sigma(\mathbf{k})$. Since $R_0$ is convex, the boundary condition is valid and it holds that for all $\mathbf{i}, \mathbf{k} \in \mathbb{Z}^d$,

$$
\begin{aligned}
|{}_sN_{\mathbf{i},n}(\mathbf{k}) - |{}_sR_n \cap \mathbf{k} + {}_sR_n|| &\leq \mathcal{C}({}_s\lambda_n^{\max})^{d-1}, \\
||{}_sR_n| - |{}_sR_n \cap \mathbf{k} + {}_sR_n|| &\leq \mathcal{C}\|\mathbf{k}\|_\infty^d({}_s\lambda_n^{\max})^{d-1}
\end{aligned}
\tag{10.12}
$$

from Lemma 9.3 and (10.9). Then (10.3), Lemma 10.4 and (10.12) give $\sum_{\mathbf{k} \in \mathbb{Z}^d} |G_n(\mathbf{k})| \leq \mathcal{C}$, $n \in \mathbb{Z}_+$; $G_n(\mathbf{k}) \to 0$ for $\mathbf{k} \in \mathbb{Z}^d$ and ${}_s\lambda_n \Psi_n = O(1)$. By the LDCT, we establish

$$\sum_{\mathbf{k} \in \mathbb{Z}^d} \frac{G_n(\mathbf{k})}{{}_s\lambda_n |R_0|} = o({}_s\lambda_n^{-1}), \qquad \mathrm{E}(\hat{\tau}_{n,\mathrm{NOL}}^2) - \tau_n^2 = \Psi_n(1 + o(1)),$$

representing the formulation of Theorem 4.2 in terms of $\Psi_n$. If $V(\mathbf{k})$ exists for each $\mathbf{k} \in \mathbb{Z}^d$, then (10.3) and (10.12) imply that we can use the LDCT again to produce

$$\Psi_n = \frac{-1}{{}_s\lambda_n |R_0|} \left( \sum_{\mathbf{k} \in \mathbb{Z}^d} V(\mathbf{k})\sigma(\mathbf{k}) \right)(1 + o(1)). \tag{10.13}$$

The proof of Theorem 4.2 for $\hat{\tau}_{n,\mathrm{NOL}}^2$ is now complete.

Consider $\hat{\tau}_{n,\mathrm{OL}}^2$. We can repeat the same steps as above to find

$$\mathrm{E}(\hat{\tau}_{n,\mathrm{OL}}^2) - \tau_n^2 = \Psi_n + \sum_{\mathbf{k} \in \mathbb{Z}^d} \frac{G_n^*(\mathbf{k})}{{}_s\lambda_n |R_0|} + o({}_s\lambda_n^{-1}), \qquad G_n^*(\mathbf{k}) = \frac{D_{\mathbf{0},n}(\mathbf{k})\sigma(\mathbf{k})}{{}_s\lambda_n^{d-1}}.$$

The same arguments for $G_n$ apply to $G_n^*$ and (10.13) remains valid when each $V(\mathbf{k})$ exists, $\mathbf{k} \in \mathbb{Z}^d$, establishing Theorem 4.2 for $\hat{\tau}_{n,\mathrm{OL}}^2$. Note as well that if $V(\mathbf{k})$ exists for each $\mathbf{k} \in \mathbb{Z}^d$, then Lemma 10.4 and Lemma 4.1 also imply the second formulation of the bias in Theorem 4.2. □

PROOF OF THEOREM 5.1. This follows from Theorems 3.1 and 4.1 and simple arguments from calculus involving minimization of a smooth function of a real variable. □

PROOF OF THEOREM 6.1. For a rectangle $T$, where $\prod_{j=1}^d (c_j, \tilde{c}_j) \subset T \subset \prod_{j=1}^d [c_j, \tilde{c}_j]$, $c_j, \tilde{c}_j \in \mathbb{R}$, define the border $\mathbf{Z}^d$-point set: $\mathcal{B}\{T\} = \bigcup_{j=1}^d \{\mathbf{s} = (s_1, \ldots, s_d)' \in \mathbf{Z}^d \cap T : s_j \in \{c_j, \tilde{c}_j\}\}$.

It holds that, for $\mathbf{k} \neq \mathbf{0}$, there exist $\mathcal{C} > 0$, $\mathrm{N}_{\mathbf{k}} \in \mathbb{Z}_+$, such that $n \geq \mathrm{N}_{\mathbf{k}}$,

$$|D_{\mathbf{i},n}(\mathbf{k})| \leq \mathcal{C}\|\mathbf{k}\|_\infty^{d-1} {}_s\lambda_n^{d-2}, \qquad \mathbf{i} \in \mathbb{Z}^d. \tag{10.14}$$

This can be shown easily by considering only volume approximations for those $\mathbf{Z}^d$ lattice point counts associated with the interior set $R_0^\circ$ [i.e., treating

OPTIMAL SPATIAL SUBSAMPLE SIZE 47$R_0^\circ$ as $R_0$ in $|D_{\mathbf{i},n}(\mathbf{k})|$] because the subtracted lattice point counts on the borders of $\tilde{R}_{\mathbf{i},n}$ and $\tilde{R}_{\mathbf{i},n} \cap \mathbf{k} + \tilde{R}_{\mathbf{i},n}$ are negligible:

$$|\mathcal{B}\{_s\lambda_n(\mathbf{i}+R_0)\}| - |\mathcal{B}\{_s\lambda_n(\mathbf{i}+R_0) \cap \mathbf{k} + {_s\lambda_n}(\mathbf{i}+R_0)\}|$$
$$\leq \mathcal{C}\|\mathbf{k}\|_\infty {_s\lambda_n}^{d-2}, \qquad \mathbf{i} \in \mathbb{Z}^d.$$

See Nordman (2002) for more details.

Applying (10.14) in place of Lemma 10.4, the same proof for Theorem 4.2 establishes Theorem 6.1. □

**Acknowledgments.** The authors thank the referees and an Associate Editor for a number of constructive comments and suggestions that significantly improved an earlier draft of this paper.

Mathematics Department
University of Wisconsin–La Crosse
La Crosse, Wisconsin 54601
USA
e-mail: nordman.dani@uwlax.edu

Department of Statistics
Iowa State University
Ames, Iowa 50011
USA
e-mail: snlahiri@iastate.edu